\documentclass[
	english
]{scrartcl}

\usepackage[T1]{fontenc}
\usepackage[utf8]{inputenc}

\usepackage{amsmath,amsfonts,amsthm,amssymb}
\usepackage[colorlinks,citecolor=blue]{hyperref}
\usepackage{doi,natbib}
\usepackage{natbib}
\usepackage{cancel}
\usepackage{relsize}

\usepackage{changes}

\usepackage{enumitem}
\setlist[enumerate]{label=(\alph*)}

\usepackage[figure]{hypcap}
\usepackage{graphicx}
\graphicspath{{./img/}}
\usepackage{subcaption} 

\usepackage{preamble}

\usepackage{marginnote}

\usepackage[capitalise,nameinlink]{cleveref}
\allowdisplaybreaks

\numberwithin{equation}{section}

\newcommand\norm[1]{\left\Vert#1\right\Vert}
\newcommand\nnorm[1]{\Vert#1\Vert}

\newcommand\dual[2]{\left\langle #1, #2\right\rangle}

\newcommand\innerprod[2]{\left\langle #1, #2\right\rangle}
\newcommand\ninnerprod[2]{\langle #1, #2\rangle}

\newcommand\N{\mathbb{N}}
\newcommand\R{\mathbb{R}}

\newcommand\tto{\rightrightarrows}

\newcommand{\dist}{\operatorname{dist}}

\newcommand{\dom}{\operatorname{dom}}
\renewcommand{\Im}{\operatorname{Im}}
\newcommand{\gph}{\operatorname{gph}}
\newcommand{\epi}{\operatorname{epi}}

\newcommand{\xb}{\bar x}

\newcommand{\yb}{\bar y}

\DeclareMathAlphabet{\mathpzc}{OT1}{pzc}{m}{it}
\newcommand\oo{\mathpzc{o}}

\newtheorem{theorem}{Theorem}[section]
\newtheorem{lemma}[theorem]{Lemma}
\newtheorem{proposition}[theorem]{Proposition}

\newtheorem{corollary}[theorem]{Corollary}
\newtheorem{remark}[theorem]{Remark}
\newtheorem{definition}[theorem]{Definition}
\newtheorem{example}[theorem]{Example}

\definecolor{mygreen}{rgb}{0.0,0.7,0.0}
\definecolor{mybrown}{rgb}{0.5,0.5,0.0}

\begin{document}

\title{On the directional asymptotic approach in optimization theory\\
	Part B: constraint qualifications}
\author{%
	Mat\'{u}\v{s} Benko%
	\footnote{%
		University of Vienna,
		Applied Mathematics and Optimization,
		1090 Vienna,
		Austria,
		\email{matus.benko@univie.ac.at},
		\url{https://www.mat.univie.ac.at/\~rabot/group.html}
		}
	\and
	Patrick Mehlitz%
	\footnote{%
		Brandenburgische Technische Universit\"at Cottbus--Senftenberg,
		Institute of Mathematics,
		03046 Cottbus,
		Germany,
		\email{mehlitz@b-tu.de},
		\url{https://www.b-tu.de/fg-optimale-steuerung/team/dr-patrick-mehlitz},
		\orcid{0000-0002-9355-850X}%
		}
	}
\publishers{}
\maketitle

\begin{abstract}
	 During the last years, asymptotic (or sequential) constraint qualifications, which postulate
	 upper semicontinuity of certain set-valued mappings and provide a natural
	 companion of asymptotic stationarity conditions, have been shown to
	 be comparatively mild, on the one hand, while possessing inherent practical relevance
	 from the viewpoint of numerical solution methods, on the other one.
	 Based on recent developments, the theory in this paper enriches asymptotic
	 constraint qualifications for very general nonsmooth optimization problems
	 over inverse images of set-valued mappings
	 by incorporating directional data. 
	 We compare these new directional asymptotic regularity conditions
	 with standard constraint qualifications from nonsmooth optimization.
	 Further, we introduce directional concepts of pseudo- and quasi-normality
	 which apply to set-valued mappings.
	 It is shown that these properties provide sufficient conditions for the validity
	 of directional asymptotic regularity.
	 Finally, a novel coderivative-like variational tool is introduced which allows 
	 to study the presence of directional asymptotic regularity.
	 For geometric constraints, it is illustrated that all appearing objects can be 
	 calculated in terms of initial problem data.
\end{abstract}


\begin{msc}	
	\mscLink{49J52}, \mscLink{49J53}, \mscLink{90C48}
\end{msc}

\section{Introduction}\label{sec:introduction}

In recent years, sequential concepts of stationarity and regularity received much attention
not only in standard nonlinear optimization,
see \cite{AndreaniMartinezSvaiter2010,AndreaniHaeserMartinez2011,AndreaniMartinezRamosSilva2016,AndreaniMartinezRamosSilva2018},
but also in complementarity-, cardinality-, and switching-constrained programming,
see \cite{AndreaniHaeserSecchinSilva2019,KanzowRaharjaSchwartz2021,LiangYe2021,Ramos2019},
conic optimization,
see \cite{AndreaniGomezHaeserMitoRamos2021},
nonsmooth optimization,
see \cite{HelouSantosSimoes2020,Mehlitz2020a,Mehlitz2022},
or even infinite-dimensional optimization,
see \cite{BoergensKanzowMehlitzWachsmuth2019,KanzowSteckWachsmuth2018,KrugerMehlitz2021}.
The interest in sequential stationarity conditions is based on the observation that they
hold at local minimizers in the absence of constraint qualifications, and that 
different types of solution algorithms like multiplier-penalty- and some SQP-methods
naturally compute such points. Sequential constraint qualifications provide conditions
which guarantee that a sequentially stationary point is already stationary in
classical sense, e.g., a Karush--Kuhn--Tucker-point 
in standard nonlinear programming or an Mordukhovich-stationary 
(M-stationary for short)
point in nonsmooth optimization. Naturally, this amounts to upper semicontinuity of
certain problem-tailored set-valued mappings. It has been reported, e.g., in
\cite{AndreaniMartinezRamosSilva2016,LiangYe2021,Mehlitz2020a,Ramos2019}
that sequential constraint qualifications are comparatively mild. 
Inherently from their construction, sequential constraint qualifications simplify
the convergence analysis of some numerical solution procedures.

In \cite{BenkoMehlitz2022}, we have shown how sequential stationarity for very general
nonsmooth problems can be enriched by directional information in terms of critical
directions at local minimizers and limiting normals in such directions, see
\cite{BenkoGfrererOutrata2019} for an introduction to as well as an overview of
the directional limiting calculus. We then used so-called pseudo-coderivatives,
see \cite{Gfrerer2014a} as well, in order to come up with mixed-order and even
M-stationarity conditions under suitable qualification conditions at local
minimizers. In the present paper, we strike a different path to benefit from
these novel findings regarding directional sequential stationarity.

The particular
sequential stationarity conditions from \cite{BenkoMehlitz2022} allow us to 
introduce directional sequential (or asymptotic) qualification conditions whose validity
directly yields M-stationarity of local minimizers,
see \cref{sec:directional_asymptotic_regularity}.  
Roughly speaking, these conditions demand certain control of unbounded input 
sequences associated with the regular coderivative of the underlying
set-valued mapping in a neighborhood of the reference point.
The directional approach
reveals that asymptotic regularity is only necessary in critical directions
and with respect to (w.r.t.) sequences satisfying some additional conditions. This way, we
can relate our new constraint qualifications with already existing ones
from the literature. Exemplary, as in \cite{Mehlitz2020a}, we observe that the concept
is independent of both, (directional) metric subregularity, see \cite{Gfrerer2013},
and the celebrated \emph{First-Order Sufficient Condition for Metric Subregularity},
see \cite{GfrererKlatte2016}.

In \cref{sec:pseudo_quasi_normality}, we introduce directional versions of
pseudo- and quasi-normality for abstract set-valued mappings. It is
illustrated that these conditions generalize the ones from 
\cite{BaiYeZhang2019,BenkoCervinkaHoheisel2019} where the authors merely
consider so-called geometric and, in particular, disjunctive constraint
systems. We show that directional pseudo- and quasi-normality are sufficient
for directional metric subregularity as well as directional asymptotic regularity.
Furthermore, we discuss how directional pseudo- and quasi-normality can be
specified for equilibrium-constrained programs which cover bilevel optimization
problems and models with (quasi-) variational inequality constraints, see e.g.\
\cite{Dempe2002,DempeKalashnikovPerezValdesKalashnykova2015,FacchneiPang2003,LuoPangRalph1996,OutrataKocvaraZowe1998}.

Finally, a new directional coderivative-like tool, the directional super-coderivative,
is introduced in \cref{sec:asymptotic_regularity_via_super_coderivative} which can
be applied beneficially for checking validity of directional asymptotic regularity.
In the presence of so-called metric pseudo-regularity, 
see \cite{Gfrerer2014a} again, this leads to conditions in terms of the aforementioned
pseudo-coderivatives. Noting that these generalized derivatives can be computed
in terms of initial problem data
for so-called feasibility mappings associated with geometric constraints,
we can specify our findings for such constraint systems.
As it turns out, the approach recovers our findings from \cite{BenkoMehlitz2022}
in different way.
Furthermore, we show that the explicit sufficient conditions for directional asymptotic regularity
provide constraint qualifications for M-stationarity which are not stronger than the
First- and Second-Order Sufficient Condition for Metric Subregularity from \cite{GfrererKlatte2016}.

The remainder of the paper is organized as follows.
In \cref{sec:notation}, we comment on the notation in this paper and recall some
fundamental tools from variational analysis and generalized differentiation which we
are going to exploit for our analysis. Furthermore, the underlying model problem from
nonsmooth optimization is introduced and the associated sequential stationarity
condition from \cite{BenkoMehlitz2022} is presented.
\Cref{sec:directional_asymptotic_regularity} is dedicated to the introduction of
directional notions of asymptotic regularity. Additionally, we comment on elementary
relations to other constraint qualifications from nonsmooth optimization. 
Some examples are used to visualize our findings.
In \cref{sec:pseudo_quasi_normality}, directional notions of pseudo- and quasi-normality
are suggested which address arbitrary set-valued mappings with a closed graph.
It is shown that these conditions serve as sufficient conditions for directional
metric subregularity and directional asymptotic regularity. Furthermore, we demonstrate
that these conditions provide suitable generalizations of already available concepts
in the literature which address geometric constraint systems. Finally, we specify
directional pseudo- and quasi-normality for constraint systems associated with
equilibrium conditions since these are modeled with the aid of set-valued mappings
in general.
Yet another way to check the presence of directional asymptotic regularity is 
presented in \cref{sec:asymptotic_regularity_via_super_coderivative}. Therein,
we define the super-coderivative of a set-valued mapping, which
is closely related to pseudo-coderivatives, and comment
on its relationship to the validity of directional asymptotic regularity. 
These findings are made precise for geometric constraint systems by
exploiting the calculus rules for the pseudo-coderivative we already obtained
in \cite{BenkoMehlitz2022}. Some concluding remarks close the paper in
\cref{sec:conclusions}.

\section{Notation and preliminaries}\label{sec:notation}

In this paper, we mainly exploit standard notation coined in
\cite{AubinFrankowska2009,BonnansShapiro2000,RockafellarWets1998,Mordukhovich2018}.

\subsection{Basic notation}\label{sec:basic_notation}

In this paper, $\mathbb X$ and $\mathbb Y$ are Euclidean spaces, i.e., finite-dimensional
Hilbert spaces, with inner product $\dual{\cdot}{\cdot}$ and corresponding norm
$\norm{\cdot}$ (the associated space will always be clear from the context).
The unit sphere of $\mathbb X$ will be denoted by $\mathbb S_{\mathbb X}$.
For a given set $Q\subset X$ and some point $\bar x$, we use
$\bar x+Q:=Q+\bar x:=\{x+\bar x\in\mathbb X\,|\,x\in Q\}$ for simplicitly.
The adjoint of a given linear operator $A\colon\mathbb X\to\mathbb Y$ will
be denoted by $A^*\colon\mathbb Y\to\mathbb X$.

For a continuously differentiable mapping $g\colon\mathbb X\to\mathbb Y$, we use
$\nabla g(\bar x)\colon\mathbb X\to\mathbb Y$ in order to denote the derivative of
$g$ at $\bar x\in\mathbb X$ which is a linear mapping between $\mathbb X$ and $\mathbb Y$.
For twice continuously differentiable $g$ and a vector $\lambda\in\mathbb Y$,
$\langle \lambda,g\rangle(x):=\langle \lambda,g(x)\rangle$ for each $x\in\mathbb X$ 
defines the associated scalarization mapping $\dual{\lambda}{g}\colon\mathbb X\to\R$.
By $\nabla\dual{\lambda}{g}(\bar x)$ and $\nabla^2\dual{\lambda}{g}(\bar x)$ we 
denote the first- and second-order derivatives of this map w.r.t.\ the variable
which enters $g$ at $\bar x$. Furthermore, for each $u\in\mathbb X$, we set
\[
	\nabla^2g(\bar x)[u,u]
	:=
	\sum_{i=1}^m\dual{u}{\nabla^2\dual{e_i}{g}(\bar x)(u)}e_i
\]
where $\{e_1,\ldots,e_m\}\subset\mathbb Y$ is the canonical basis of $\mathbb Y$.

\subsection{Variational analysis and generalized differentiation}\label{sec:generalized_differentiation}

Let us fix a closed set $Q\subset\mathbb X$ and some point $\bar x\in Q$.
The set
\[
	\mathcal T_Q(\bar x)
	:=
	\left\{
		d\in\mathbb X\,\middle|\,
			\begin{aligned}
				&\exists\{d_k\}_{k\in\N}\subset\mathbb X,\,\exists\{t_k\}_{k\in\N}\subset\R_+
				\colon\\
				&\qquad
				d_k\to d,\,t_k\searrow 0,\,\bar x+t_kd_k\in Q\,\forall k\in\N
			\end{aligned}
	\right\}
\]
is called the \emph{tangent} (or Bouligand) cone to $Q$ at $\bar x$.
Furthermore, we exploit
\begin{align*}
	\widehat{\mathcal N}_Q(\bar x)
	&:=
	\left\{
		\eta\in\mathbb X\,\middle|\,
			\forall x\in Q\colon\,\dual{\eta}{x-\bar x}\leq\oo(\norm{x-\bar x})
	\right\},
	\\
	\mathcal N_Q(\bar x)
	&:=
	\left\{
		\eta\in\mathbb X\,\middle|\,
			\begin{aligned}
				&\exists\{x_k\}_{k\in\N}\subset Q,\,\exists\{\eta_k\}_{k\in\N}\subset\mathbb X
				\colon
				\\
				&\qquad x_k\to\bar x,\,\eta_k\to\eta,\,
				\eta_k\in\widehat{\mathcal N}_Q(x_k)\,\forall k\in\N
			\end{aligned}
	\right\}
\end{align*}
which are referred to as the \emph{regular} (or Fr\'{e}chet) and
\emph{limiting} (or Mordukhovich) \emph{normal cone} to $Q$ at $\bar x$ in the literature.
It is well known that these cones coincide with the normal cone
in the sense of convex analysis whenever $Q$ is a convex set.
For the purpose of completeness, for each $\tilde x\notin Q$, we put
$\mathcal T_Q(\tilde x):=\varnothing$ and
$\widehat{\mathcal N}_Q(\tilde x)=\mathcal N_Q(\tilde x):=\varnothing$.

For some direction $u\in\mathbb X$, we make use of
\[
	\mathcal N_Q(\bar x;u)
	:=
	\left\{
		\eta\in\mathbb X\,\middle|\,
		\begin{aligned}
			&\exists\{u_k\}_{k\in\N}\subset\mathbb X,\,
			\exists\{t_k\}_{k\in\N}\subset\R_+,\,
			\exists\{\eta_k\}_{k\in\N}\subset\mathbb X\colon\\
			&\qquad u_k\to u,\,t_k\searrow 0,\,\eta_k\to\eta,\,
			\eta_k\in\widehat{\mathcal N}_Q(\bar x+t_ku_k)\,\forall k\in\N
		\end{aligned}
	\right\}
\]
which is called the \emph{limiting normal cone} to $Q$ at $\bar x$
\emph{in direction $u$}. Note that this set is empty when $\bar x\notin Q$ or
$u\notin\mathcal T_Q(\bar x)$.
In case where $Q$ is convex, we obtain 
$\mathcal N_Q(\bar x;u)=\mathcal N_Q(\bar x)\cap[u]^\perp$
where $[u]^\perp:=\{\eta\in\mathbb X\,|\,\innerprod{\eta}{u}=0\}$
is the annihilator of $u$.

The limiting normal cone to a set is well known for its robustness, i.e., 
it is outer semicontinuous as a set-valued mapping.
In the course of the paper, we exploit an analogous property of the directional limiting 
normal cone which has been validated in \cite[Proposition~2]{GfrererYeZhou2022}.
\begin{lemma}\label{lem:robustness_directional_limiting_normals}
	Let $Q\subset\mathbb X$ be closed and fix $\bar x\in Q$.
	Then, for each $u\in\mathbb X$, we have
	\[
	\mathcal N_Q(\bar x;u)
	=
	\left\{
		\eta\in\mathbb X\,\middle|\,
			\begin{aligned}
				&\exists\{u_k\}_{k\in\N}\subset\mathbb X,\,
					\exists\{t_k\}_{k\in\N}\subset\R_+,\,
					\exists\{\eta_k\}_{k\in\N}\subset\mathbb X\colon\\
				&\qquad u_k\to u,\,t_k\searrow 0,\,\eta_k\to\eta,\,
				\eta_k\in \mathcal N_Q(\bar x+t_ku_k)\,\forall k\in\N
			\end{aligned}
	\right\}.
\]
\end{lemma}

Next, we recall some fundamental notions of generalized differentiation.
Let us start with a locally Lipschitz continuous function $\varphi\colon\mathbb X\to\R$
and fix $\bar x\in\mathbb X$. The sets
\begin{align*}
	\widehat\partial\varphi(\bar x)
	&:=
	\left\{
		\eta\in\mathbb X\,\middle|\,
			(\eta,-1)\in\widehat{\mathcal N}_{\epi\varphi}(\bar x,\varphi(\bar x))
	\right\},
	\\
	\partial\varphi(\bar x)
	&:=
	\left\{
		\eta\in\mathbb X\,\middle|\,
			(\eta,-1)\in\mathcal N_{\epi\varphi}(\bar x,\varphi(\bar x))
	\right\}
\end{align*}
are referred to as the \emph{regular} and \emph{limiting subdifferential}
of $\varphi$ at $\bar x$. 
Here, $\epi\varphi:=\{(x,\alpha)\in\mathbb X\times\R\,|\,\varphi(x)\leq\alpha\}$
denotes the epigraph of $\varphi$.
Furthermore, for some direction $u\in\mathbb X$,
\[
	\partial\varphi(\bar x;u)
	:=
	\left\{
		\eta\in\mathbb X\,\middle|\,
			\begin{aligned}
				&\exists\{u_k\}_{k\in\N}\subset\mathbb X,\,
				\exists\{t_k\}_{k\in\N}\subset\R_+,\,
				\exists\{\eta_k\}_{k\in\N}\subset\mathbb X\colon
				\\
				&\qquad
				u_k\to u,\,t_k\searrow 0,\,\eta_k\to\eta,\,
				\eta_k\in\widehat\partial\varphi(\bar x+t_ku_k)\,\forall k\in\N
			\end{aligned}
	\right\}
\]
is referred to as the \emph{limiting subdifferential} of $\varphi$ at $\bar x$
\emph{in direction $u$}.

Let $\Phi\colon\mathbb X\tto\mathbb Y$ be a set-valued mapping.
The sets 
$\dom\Phi:=\{x\in\mathbb X\,|\,\Phi(x)\neq\varnothing\}$,
$\gph\Phi:=\{(x,y)\in\mathbb X\times\mathbb Y\,|\,y\in\Phi(x)\}$,
$\ker\Phi:=\{x\in\mathbb X\,|\,0\in\Phi(x)\}$,
and
$\Im\Phi:=\bigcup_{x\in\mathbb X}\Phi(x)$
are called the domain, graph, kernel, and image of $\Phi$, respectively.

We fix some point $(\bar x,\bar y)\in\gph\Phi$.
The set-valued mapping $D\Phi(\bar x,\bar y)\colon\mathbb X\tto\mathbb Y$ given by
$\gph D\Phi(\bar x,\bar y):=\mathcal T_{\gph\Phi}(\bar x,\bar y)$ is called
\emph{graphical derivative} of $\Phi$ at $(\bar x,\bar y)$.
In case where $\Phi$ is single-valued at $\bar x$, we exploit
$D\Phi(\bar x)\colon\mathbb X\tto\mathbb Y$ for brevity
of notation.
In \cite[Definition~2.4]{BenkoMehlitz2022}, we introduced the so-called
\emph{graphical subderivative} of $\Phi$ at $(\bar x,\bar y)$ to be the 
set-valued mapping 
$D_\textup{sub}\Phi(\bar x,\bar y)\colon\mathbb S_{\mathbb X}\tto\mathbb S_{\mathbb Y}$
which assigns to every $u\in\mathbb S_{\mathbb X}$ the set of vectors 
$v\in\mathbb S_{\mathbb Y}$
such that there are sequences $\{u_k\}_{k\in\N}\subset\mathbb X$,
$\{v_k\}_{k\in\N}\subset\mathbb Y$, and $\{t_k\}_{k\in\N},\{\tau_k\}_{k\in\N}\subset\R_+$
such that $u_k\to u$, $v_k\to v$, $t_k\searrow 0$, $\tau_k\searrow 0$,
$\tau_k/t_k\to\infty$, and $(\bar x+t_ku_k,\bar y+\tau_kv_k)\in\gph\Phi$ for all $k\in\N$.
Some fundamental calculus rules for the graphical subderivative, particularly,
regarding so-called normal cone mappings, can be found in \cite[Section~2.3]{BenkoMehlitz2022}.

Let us now turn our attention to dual concepts of generalized differentiation which
will be of essential importance in this paper.
We refer to $\widehat D^*\Phi(\bar x,\bar y),D^*\Phi(\bar x,\bar y)\colon\mathbb Y\tto\mathbb X$
given by
\begin{align*}
	\widehat D^*\Phi(\bar x,\bar y)(y^*)
	&:=
	\left\{
		x^*\in\mathbb X\,\middle|\,(x^*,-y^*)\in\widehat{\mathcal N}_{\gph\Phi}(\bar x,\bar y)
	\right\},
	\\
	D^*\Phi(\bar x,\bar y)(y^*)
	&:=
	\left\{
		x^*\in\mathbb X\,\middle|\,(x^*,-y^*)\in\mathcal N_{\gph\Phi}(\bar x,\bar y)
	\right\}
\end{align*}
for each $y^*\in\mathbb Y$ as \emph{regular} and \emph{limiting coderivative} of
$\Phi$ at $(\bar x,\bar y)$. 
For a pair of directions $(u,v)\in\mathbb X\times\mathbb Y$, the set-valued mapping
$D^*\Phi((\bar x,\bar y);(u,v))\colon\mathbb Y\tto\mathbb X$ defined by means of
\begin{align*}
	D^*\Phi((\bar x,\bar y);(u,v))(y^*)
	:=
	\left\{x^*\in\mathbb X\,\middle|\,
		\begin{aligned}
			&\exists\{u_k\}_{k\in\N},\{x_k^*\}_{k\in\N}\subset\mathbb X,\,
			\exists\{v_k\}_{k\in\N},\{y_k^*\}_{k\in\N}\subset\mathbb Y,\\
			&\exists\{t_k\}_{k\in\N}\subset\R_+\colon
			\\
			&\qquad u_k\to u,\,v_k\to v,\,t_k\searrow 0,\,x_k^*\to x^*,\,y_k^*\to y^*,\\
			&\qquad
				x_k^*\in\widehat D^*\Phi(\bar x+t_ku_k,\bar y+t_kv_k)(y_k^*)\,\forall k\in\N
		\end{aligned}
	\right\}
\end{align*}
for all $y^*\in\mathbb Y$ is called the \emph{limiting coderivative} of $\Phi$
at $(\bar x,\bar y)$ \emph{in direction $(u,v)$}.
We note that the latter is only reasonable if $v\in D\Phi(\bar x,\bar y)(u)$, and that
\[
	D^*\Phi((\bar x,\bar y);(u,v))(y^*)
	=
	\left\{
		x^*\in\mathbb X\,\middle|\,(x^*,-y^*)\in\mathcal N_{\gph\Phi}((\bar x,\bar y);(u,v))
	\right\}
\]
holds for all $y^*\in\mathbb Y$.
Again, we exploit
$\widehat D^*\Phi(\bar x),D^*\Phi(\bar x),D^*\Phi(\bar x;(u,v))\colon\mathbb Y\tto\mathbb X$
for brevity if $\Phi$ is single-valued at $\bar x$.

For a given order $\gamma>1$ and $(u,v)\in\mathbb S_{\mathbb X}\times\mathbb Y$, 
let the \emph{pseudo-coderivative} of $\Phi$ 
\emph{of order $\gamma$} at $(\bar x,\bar y)$ \emph{in direction $(u,v)$}
be the set-valued mapping $D^*_\gamma\Phi((\bar x,\bar y);(u,v))\colon\mathbb Y\tto\mathbb X$
which assigns to $y^*\in\mathbb Y$ the set of all vectors $x^*\in\mathbb X$ such that
there exist sequences $\{u_k\}_{k\in\N},\{x_k^*\}_{k\in\N}\subset\mathbb X$,
$\{v_k\}_{k\in\N},\{y_k^*\}_{k\in\N}\subset\mathbb Y$, and $\{t_k\}_{k\in\N}\subset\R_+$
such that $u_k\to u$, $v_k\to v$, $x_k^*\to x^*$, $y_k^*\to y^*$, $t_k\searrow 0$,
and
\[
	\forall k\in\N\colon\quad
	(t_k\norm{u_k})^{\gamma-1}x_k^*
	\in 
	\widehat{D}^*\Phi(\bar x+t_ku_k,\bar y+(t_k\norm{u_k})^\gamma v_k)(y_k^*).
\]
This notion of a pseudo-coderivative originates from \cite{BenkoMehlitz2022} where it has
been used to derive mixed-order stationarity conditions for nonsmooth optimization
problems. However, the concept of pseudo-coderivative is a little older and dates back
to \cite{Gfrerer2014a}. Therein, a set-valued mapping 
$\widetilde D^*_\gamma((\bar x,\bar y);(u,v))\colon\mathbb Y\tto\mathbb X$
is called pseudo-coderivative of $\Phi$ of order $\gamma$ at $(\bar x,\bar y)$ in
direction $(u,v)$ when there exist sequences
$\{u_k\}_{k\in\N},\{x_k^*\}_{k\in\N}\subset\mathbb X$,
$\{v_k\}_{k\in\N},\{y_k^*\}_{k\in\N}\subset\mathbb Y$, and $\{t_k\}_{k\in\N}\subset\R_+$
such that $u_k\to u$, $v_k\to v$, $x_k^*\to x^*$, $y_k^*\to y^*$, $t_k\searrow 0$,
and
\[
	\forall k\in\N\colon\quad
	(t_k\norm{u_k})^{\gamma-1}x_k^*
	\in 
	\widehat{D}^*\Phi(\bar x+t_ku_k,\bar y+t_kv_k)(y_k^*).
\]
In order to distinguish both concepts, we use the slightly different notation from above
and refer to $\widetilde D^*_\gamma((\bar x,\bar y);(u,v))$ as 
\emph{Gfrerer's directional pseudo-coderivative}.
By definition of these tools, we have the trivial estimate
\[
	D^*_\gamma\Phi((\bar x,\bar y);(u,v))(y^*)
	\subset
	\widetilde D^*_\gamma\Phi((\bar x,\bar y);(u,0))(y^*)
\]
for all $y^*\in\mathbb Y$.

Let us recall that $\Phi$ is said to be \emph{metrically subregular}
at $(\bar x,\bar y)$ in direction $u\in\mathbb X$ whenever there
are constants $\varepsilon>0$, $\delta>0$, and $\kappa>0$ such that
\[
	\forall x\in\bar x+\mathbb B_{\varepsilon,\delta}(u)\colon\quad
	\dist(x,\Phi^{-1}(\bar y))
	\leq
	\kappa\,\dist(\bar y,\Phi(x)),
\]
where 
$\mathbb B_{\varepsilon,\delta}(u)
:=\{v\in\mathbb X\,|\,\norm{\norm{v}u-\norm{u}v}\leq\delta\norm{u}\norm{v},\,
\norm{v}\leq\varepsilon\}$
is a so-called \emph{directional neighborhood} of $u$
and 
$\Phi^{-1}(\bar y):=\{x\in\mathbb X\,|\,\bar y\in\Phi(x)\}$
is the inverse image of $\bar y$ under $\Phi$.
In case where this is fulfilled for $u:=0$, $\Phi$ is said to be
metrically subregular at $(\bar x,\bar y)$.

Coderivatives have turned out to be suitable tools in order to characterize
local Lipschitz or regularity properties of set-valued mappings.
Exemplary, let us mention that the so-called \emph{Mordukhovich criterion}
\[
	\ker D^*\Phi(\bar x,\bar y)=\{0\}
\]
is equivalent to $\Phi$ being \emph{metrically regular} at $(\bar x,\bar y)$,
see e.g.\ \cite[Section~3.1]{Mordukhovich2018} for a definition and this result.
Furthermore, the condition
\[
	\forall u\in\mathbb S_{\mathbb X}\colon\quad
	\ker D^*\Phi((\bar x,\bar y);(u,0))=\{0\}
\]
is sufficient for $\Phi$ to be metrically subregular at $(\bar x,\bar y)$,
see \cite{Gfrerer2013}, which is why it is called 
\emph{First-Order Sufficient Condition For Metric Subregularity}
(FOSCMS for short)
in the literature. We also note that, for some fixed $u\in\mathbb S_{\mathbb X}$,
$\ker D^*\Phi((\bar x,\bar y);(u,0))=\{0\}$
implies that $\Phi$ is metrically subregular at $(\bar x,\bar y)$
in direction $u$.
An analogous sufficient condition for \emph{metric pseudo-subregularity} of order $\gamma > 1$
in terms of Gfrerer's directional pseudo-coderivative was derived in \cite{Gfrerer2014a},
and in \cite[Lemma~2.9]{BenkoMehlitz2022} it was slightly modified to
\[
	\forall u\in\mathbb S_{\mathbb X}\colon\quad
	\ker D^*_\gamma\Phi((\bar x,\bar y);(u,0))=\{0\}.
\]

Finally, we would like to provide some basic calculus rules for the coderivative
of so-called constraint mappings.
Therefore, recall that some single-valued function $g\colon\mathbb X\to\mathbb Y$ is called \emph{calm} at $x\in\mathbb X$ whenever there 
are a neighborhood $U\subset\mathbb X$ of $x$ and a constant $L>0$ such that
\[
	\forall x'\in U\colon\quad \nnorm{g(x')-g(x)}\leq L\nnorm{x'-x}.
\]
Furthermore, for some direction $u$, $g$ is referred to as \emph{calm in direction $u$}
at $x$ if there are constants $\varepsilon>0$, $\delta>0$, and $L>0$ such that
\[
	\forall x'\in\bar x+\mathbb B_{\varepsilon,\delta}(u)\colon\quad
	\nnorm{g(x')-g(x)}\leq L\nnorm{x'-x}.
\]

\begin{lemma}\label{lem:coderivatives_constraint_maps}
	Let $g\colon\mathbb X\to\mathbb Y$ be continuous, and let $D\subset\mathbb Y$
	be nonempty as well as closed.
	We consider the constraint map $\Phi\colon\mathbb X\tto\mathbb Y$
	given by $\Phi(x):=g(x)-D$ for all $x\in\mathbb X$.
	Fix $(x,y)\in\gph\Phi$. Then the following statements hold.
	\begin{enumerate}
		\item\label{item:constraint_maps_regular_coderivative} 
		For each $y^*\in\mathbb Y$, we have
        \[
		\widehat{D}^*\Phi(x,y)(y^*)
		\subset
		\begin{cases}
			\widehat{D}^*g(x) (y^*)	&	y^*\in \widehat{\mathcal N}_D(g(x)-y),\\
			\varnothing				&	\text{otherwise,}
		\end{cases}
		\]
        and the opposite inclusion holds if $g$ is calm at $x$.
		\item\label{item:constraint_maps_limiting_coderivative}
		 For each $y^*\in\mathbb Y$, we have
		\[
		D^*\Phi(x,y)(y^*)
		\subset
		\begin{cases}
			D^*g(x) (y^*)	&	y^*\in \mathcal N_D(g(x)-y),\\
			\varnothing	&	\text{otherwise}.
		\end{cases}
		\]
		\item\label{item:constraint_maps_directional_limiting_coderivative} 
		 For each pair of directions $(u,v)\in\mathbb X\times\mathbb Y$
			and each $y^*\in\mathbb Y$, we have
			\[
				D^*\Phi((x,y);(u,v))(y^*)
				\subset
				\begin{cases}
					\bigcup\limits_{w\in Dg(x)(u)}
					D^*g(x;(u,w))(y^*)	& y^*\in\mathcal N_D(g(x)-y;w-v),\\
					\varnothing				&\text{otherwise}
				\end{cases}
			\]
			provided $g$ is calm at $x$.
	\end{enumerate}
\end{lemma}
\begin{proof}
	\begin{enumerate}
		\item For the proof of the statement, we observe that
			$\gph\Phi=\gph g-(\{0\}\times D)$ is valid.
			Now, we exploit the sum rule from \cite{BenkoMehlitz2020}.
			Therefore, let us introduce the surrogate mapping 
			$M\colon\mathbb X\times\mathbb Y\tto(\mathbb X\times\mathbb Y)\times(\mathbb X\times\mathbb Y)$ given by
			\begin{equation}\label{eq:intermediate_map_sum_rule}
				\begin{aligned}
				M(x,y)
				&:=
				\{((\tilde x,g(\tilde x)),(0,\tilde y))\,|\,\tilde x=x,\,\tilde y\in -D, y=g(\tilde x)+\tilde y\} 
				\\
				&=
				\begin{cases}
					\{((x,g(x)),(0,y-g(x)))\}	& g(x)-y\in D,\\
					\varnothing			& \text{otherwise}
				\end{cases}
				\end{aligned}
			\end{equation}
			for all $(x,y)\in\mathbb X\times\mathbb Y$,
			and observe that $\gph\Phi=\dom M$ holds while $M$ is single-valued and continuous on $\gph\Phi$.
			Now, we find
			\[
				\widehat{\mathcal N}_{\gph\Phi}(x,y)
				\subset
				\widehat{D}^*M((x,y),((x,g(x)),(0,y-g(x))))((0,0),(0,0))
			\]
			for all $(x,y)\in\gph\Phi$
			from \cite[Theorem~3.1]{BenkoMehlitz2020}, and the converse
			inclusion holds if $g$ is calm at $x$ since this ensures that $M$ is
			so-called isolatedly calm at the point of interest, 
			see \cite[Corollary~4.4, Section~5.1.1]{BenkoMehlitz2020}.
			Now, computing the regular coderivative of $M$ via 
			\cite[Lemmas~2.1, 2.2]{BenkoMehlitz2020} yields the claim.
		\item
			The proof is similar as the one of the first statement.
			Again, we exploit the mapping $M$ given in \eqref{eq:intermediate_map_sum_rule} and
			apply \cite[Theorem~3.1]{BenkoMehlitz2020} while observing that $M$ is 
			so-called inner semicompact w.r.t.\ its domain 
			at each point $(x,y)\in\gph\Phi$ by continuity of $g$.
		\item 
			This assertion can be shown in similar way as the second one.
	\end{enumerate}
\end{proof}

Let us note that the upper estimate in~\ref{item:constraint_maps_regular_coderivative} was also shown in 
\cite[Lemma~3.2]{BaiYeZhang2019}, but it actually follows directly from
\cite[Exercise~6.44]{RockafellarWets1998} upon realizing $\gph\Phi=\gph g-(\{0\}\times D)$.
In case where $g$ is not calm at the reference point, one can still obtain an upper
estimate for the directional limiting coderivative from \cite[Theorem~3.1]{BenkoMehlitz2020}
which is slightly more technical since it comprises another union over
$w\in Dg(x)(0)\cap\mathbb S_{\mathbb Y}$.

\subsection{The model problem}\label{sec:model_program}

Let $\varphi\colon\mathbb X\to\R$ be a locally Lipschitz continuous mapping, 
assume that $\Phi\colon\mathbb X\tto\mathbb Y$ has a closed graph,
and fix $\bar y\in\Im\Phi$.
In this paper,
we investigate the rather general nonsmooth optimization problem
\begin{equation}\label{eq:nonsmooth_problem}\tag{P}
	\min\{\varphi(x)\,|\,\bar y\in\Phi(x)\}.
\end{equation}
The feasible set of \eqref{eq:nonsmooth_problem} will be denoted by 
$\mathcal F\subset\mathbb X$ and is, by $\bar y\in\Im\Phi$, nonempty.
Let us remark that the model \eqref{eq:nonsmooth_problem} 
covers diverse classes of optimization problems from the literature including
standard nonlinear problems, problems with geometric (particularly, disjunctive or conic) constraints, problems with (quasi-) variational inequality constraints, and
bilevel optimization problems. 
Optimality conditions and constraint qualifications for problems of this type can be
found, e.g., in \cite{Gfrerer2013,Mehlitz2020a,Mordukhovich2006,YeYe1997}.
A standard notion of stationarity, which applies to \eqref{eq:nonsmooth_problem} and
is based on the tools of limiting variational analysis, is the
one of M-stationarity.

\begin{definition}\label{def:M_stationarity}
	A feasible point $\bar x\in\mathcal F$ of \eqref{eq:nonsmooth_problem} 
	is called \emph{M-stationary} whenever there is a multiplier $\lambda\in\mathbb Y$ such that
	\[
		0\in\partial\varphi(\bar x)+D^*\Phi(\bar x,\bar y)(\lambda).
	\]
\end{definition}

In the following lemma, whose proof is analogous to the one of
\cite[Lemma~3.1]{BaiYe2021}, we point out that directional metric subregularity
of $\Phi$ implies that penalizing the constraint in \eqref{eq:nonsmooth_problem}
with the aid of the distance function yields a \emph{directionally} exact penalty function.
\begin{lemma}\label{lem:directional_exact_penalization}
	Let $\bar x\in\mathcal F$ be a local minimizer of \eqref{eq:nonsmooth_problem},
	and assume that $\Phi$ is metrically subregular at $(\bar x,\bar y)$
	in direction $u\in\mathbb S_{\mathbb X}$. 
	Then there are constants $\varepsilon>0$, $\delta>0$, and $C>0$ such
	that $\bar x$ is a local minimizer of
	\begin{equation}\label{eq:directionally_penalized_problem}
		\min\{
			\varphi(x)+C\dist(\bar y,\Phi(x))
			\,|\,
			x\in\bar x+\mathbb B_{\varepsilon,\delta}(u)
		\}.
	\end{equation}
\end{lemma} 

Let us note that this result refines well-known theory about classical 
exact penalization in the presence of metric subregularity, see e.g.\ \cite{Burke1991,Clarke1983,KlatteKummer2002}.

In order to state one of the essential findings of \cite{BenkoMehlitz2022} which provides
the basis of our investigations, we need to recall the notion of critical directions 
associated with \eqref{eq:nonsmooth_problem}.
\begin{definition}\label{def:critical_direction}
	For some feasible point $\bar x\in\mathcal F$ of \eqref{eq:nonsmooth_problem}, 
	a direction $u\in\mathbb S_{\mathbb X}$ is called \emph{critical}
	for \eqref{eq:nonsmooth_problem} at $\bar x$ 
	whenever there are sequences $\{u_k\}_{k\in\N}\subset\mathbb X$, 
	$\{v_k\}_{k\in\N}\subset\mathbb Y$,
	and $\{t_k\}_{k\in\N}\subset\R_+$ such that 
	$u_k\to u$, $v_k\to 0$, $t_k\searrow 0$, 
	and $(\bar x+t_ku_k,\bar y+t_kv_k)\in\gph\Phi$ for all $k\in\N$ as well as
	\[
		\limsup\limits_{k\to\infty}\frac{\varphi(\bar x+t_ku_k)-\varphi(\bar x)}{t_k}\leq 0.
	\]
\end{definition}

Let us note that whenever $\varphi$ is directionally differentiable 
at $\bar x\in\mathcal F$, then
$u\in\mathbb S_{\mathbb X}$ is a critical for \eqref{eq:nonsmooth_problem} 
at $\bar x$ if and only if
$\varphi'(\bar x;u)\leq 0$ and $0\in D\Phi(\bar x,\bar y)(u)$.

A directionally refined concept of M-stationarity has been shown to
serve as a necessary optimality condition under validity of directional
metric subregularity in \cite[Theorem~7]{Gfrerer2013}.
\begin{lemma}\label{lem:directional_M_stationarity_via_metric_subregularity}
	Let $\bar x\in\mathcal F$ be a local minimizer of \eqref{eq:nonsmooth_problem},
	let $u\in\mathbb S_{\mathbb X}$ be a critical direction for \eqref{eq:nonsmooth_problem}
	at $\bar x$, and let $\Phi$ be metrically subregular at $(\bar x,\bar y)$ in
	direction $u$. Then there is a multiplier $\lambda\in\mathbb Y$ such that
	\[
		0\in\partial\varphi(\bar x;u)+D^*\Phi((\bar x,\bar y);(u,0))(\lambda).
	\]
	Particularly, $\bar x$ is M-stationary.
\end{lemma}

Let us note that the above result can also be distilled 
from \cref{lem:directional_exact_penalization} by
following ideas used to prove \cite[Theorem~3.1]{BaiYe2021}.

The following result is taken from \cite[Corollary~4.4]{BenkoMehlitz2022}
and sharpens the information provided by \cite[Theorem~3.2]{Mehlitz2020a}
or \cite[Theorem~4.1]{KrugerMehlitz2021}.

\begin{theorem}\label{thm:directional_asymptotic_stationarity}
	Let $\bar x\in\mathcal F$ be a local minimizer of \eqref{eq:nonsmooth_problem}.
	Then $\bar x$ is M-stationary or there exist a critical direction 
	$u\in\mathbb S_{\mathbb X}$ for \eqref{eq:nonsmooth_problem} at $\bar x$
	and $y^* \in \mathbb Y$ 
	as well as sequences
	$\{x_k\}_{k\in\N},\{x_k'\}_{k\in\N},\{\eta_k\}_{k\in\N}\subset\mathbb X$ and
	$\{y_k\}_{k\in\N},\{y_k^*\}_{k\in\N}\subset\mathbb Y$ such that
	$x_k,x_k'\notin\Phi^{-1}(\bar y)$, $y_k\neq\bar y$, and $y_k^*\neq 0$ for all $k\in\N$,
	\begin{subequations}\label{eq:convergences_gamma=1}
		\begin{align}
			\label{eq:convergences_gamma=1_basic}
				x_k, x_k'&\to\bar x,&		\qquad	y_k&\to\bar y,&	\qquad	\eta_k&\to 0,&
				\\
			\label{eq:convergences_gamma=1_directional}
				\frac{x_k-\bar x}{\nnorm{x_k-\bar x}}&\to u,&	\qquad
				\frac{x_k'-\bar x}{\nnorm{x_k'-\bar x}}&\to u,&	\qquad
				\frac{y_k-\bar y}{\nnorm{x_k-\bar x}}&\to 0,&
				\\
			\label{eq:convergences_gamma=1_multiplier}
			y_k^*&\to y^*,&\qquad
			\frac{\norm{x_k - \bar x}}{\norm{y_k - \bar y}}\norm{y_k^*} &\to \infty,&\qquad
			\frac{y_k-\bar y}{\nnorm{y_k-\bar y}}-\frac{y_k^*}{\nnorm{y_k^*}}&\to 0,&
		\end{align}
	\end{subequations}	
	and
	\begin{equation}\label{eq:asymptotic_stationarity_gamma=1}
		\forall k\in\N\colon\quad
		\eta_k\in\widehat{\partial}\varphi(x_k')+\widehat{D}^*\Phi(x_k,y_k)\left(\frac{\norm{x_k - \bar x}}{\norm{y_k - \bar y}}y_k^*\right).
	\end{equation}
\end{theorem}

Observe that \cref{thm:directional_asymptotic_stationarity} yields a necessary optimality
condition for \eqref{eq:nonsmooth_problem} which holds in the absence of any constraint
qualification. Basically, \cref{thm:directional_asymptotic_stationarity} says that
either a local minimizer of \eqref{eq:nonsmooth_problem} is M-stationary or the
so-called approximate (or asymptotic) stationarity condition
\eqref{eq:asymptotic_stationarity_gamma=1} holds along certain sequences such that
the involved sequence of multiplier estimates given by
\[
	\forall k\in\N\colon\quad
	\lambda_k:=y_k^*\norm{x_k-\bar x}/\norm{y_k-\bar y}
\]
is unbounded.
Note that in case where $\{\lambda_k\}_{k\in\N}$ would be bounded, one could simply
take the limit in \eqref{eq:asymptotic_stationarity_gamma=1} along a suitable
subsequence and, respecting the
convergences from \eqref{eq:convergences_gamma=1}, would end up with M-stationarity
again. Thus, divergence of the multiplier estimates is natural in 
\cref{thm:directional_asymptotic_stationarity} since not all local minimizers of
\eqref{eq:nonsmooth_problem} are M-stationary in general,
see \cite[Lemma~3.4]{Mehlitz2020a} as well.

The sequential information from \eqref{eq:convergences_gamma=1} describes in great detail what must ``go wrong'' if M-stationarity fails.
We will refer to \eqref{eq:convergences_gamma=1_basic}-\eqref{eq:convergences_gamma=1_multiplier} as
basic, directional, and multiplier (sequential) information, respectively.
Clearly, one can secure M-stationarity of a local minimizer by ruling out the second alternative in \cref{thm:directional_asymptotic_stationarity} and,
as we will show, various known constraint qualifications for M-stationarity indeed do precisely that.
Let us mention here two such conditions.
Rescaling \eqref{eq:asymptotic_stationarity_gamma=1} by
$\norm{\lambda_k}$ and taking the limit $k\to\infty$ leads to a contradiction with the Mordukhovich criterion/metric regularity of $\Phi$ at $(\bar x,\bar y)$.
Respecting also the directional information \eqref{eq:convergences_gamma=1_directional} yields a contradiction with FOSCMS at $(\bar x,\bar y)$.
Thus, we obtain a result related to \cref{lem:directional_M_stationarity_via_metric_subregularity}, see \cite[Theorem~4.3]{BenkoMehlitz2022} as well.
The advantage of these two conditions lies in their simplicity, since they can be expressed via suitable derivatives, but they are a bit more restrictive.

In both cases, we have essentially discarded the multiplier information \eqref{eq:convergences_gamma=1_multiplier}
which deserves some remarks.
We have used $\norm{\lambda_k} \to \infty$, but this information is not really very important, since as we already explained,
if the multipliers remain bounded, we end up with M-stationarity anyway.
The fact that $\{y_k^*\}_{k\in\N}$ converges tells us how fast the multipliers 
$\{\lambda_k\}_{k\in\N}$ grow since we have 
$y_k^* = \lambda_k \norm{y_k-\bar y}/\norm{x_k-\bar x}$ for each $k\in\N$.
In \cref{sec:asymptotic_regularity_via_super_coderivative}, we introduce the so-called super-coderivative which incorporates this information.

Finally, $(y_k-\bar y)/\nnorm{y_k-\bar y} - \lambda_k/\nnorm{\lambda_k} \to 0$,
which equals
$(y_k-\bar y)/\nnorm{y_k-\bar y} - y_k^*/\nnorm{y_k^*} \to 0$,
means that the multipliers
precisely capture the direction from which $\{y_k\}_{k\in\N}$ converges to $\bar y$.
Equivalently, it can be expressed via $\dual{\lambda_k/\nnorm{\lambda_k}}{(y_k-\bar y)/\nnorm{y_k-\bar y}} \to 1$,
which was used in the sufficient condition for metric subregularity in \cite[Corollary 1]{Gfrerer2014a}.
This information is behind the notions of pseudo- and quasi-normality and we discuss it in detail in \cref{sec:pseudo_quasi_normality}
and also utilize it in \cref{sec:asymptotic_regularity_via_super_coderivative} to some extent.

Let us also mention that if the regular subdifferential and coderivative are replaced by the limiting ones in \eqref{eq:asymptotic_stationarity_gamma=1},
$y_k$ and $y_k^*$ can be chosen such that $(y_k-\bar y)/\nnorm{y_k-\bar y} = y_k^*/\nnorm{y_k^*}$ holds for all $k\in\N$.
This owes to the fact that the fuzzy calculus used in the proof of \cite[Theorem~4.3]{BenkoMehlitz2022}
can be replaced by the exact calculus for limiting subdifferentials.

\section{Directional asymptotic regularity in nonsmooth optimization}\label{sec:directional_asymptotic_regularity}

Based on \cref{thm:directional_asymptotic_stationarity}, the following definition introduces
concepts which may serve as (directional) qualification conditions for \eqref{eq:nonsmooth_problem}.
\begin{definition}\label{def:asymptotic_regularity}
	Let $(\bar x,\bar y)\in\gph\Phi$ be fixed.
	\begin{enumerate}
		\item The map $\Phi$ is said to be \emph{asymptotically regular at $(\bar x,\bar y)$}
			whenever the following condition holds:
			for every sequences $\{(x_k,y_k)\}_{k\in\N}\subset\gph\Phi$,
			$\{x_k^*\}_{k\in\N}\subset\mathbb X$, and $\{\lambda_k\}_{k\in\N}\subset\mathbb Y$
			as well as $x^*\in\mathbb X$ satisfying $x_k\to\bar x$, $y_k\to\bar y$,
			$x_k^*\to x^*$, and $x_k^*\in\widehat{D}^*\Phi(x_k,y_k)(\lambda_k)$ for all
			$k\in\N$, we find $x^*\in\Im D^*\Phi(\bar x,\bar y)$.
		\item\label{item:def_dir_asymp_reg} 
			For the fixed direction $u\in\mathbb S_{\mathbb X}$, 
			$\Phi$ is said to be \emph{asymptotically regular 
			at $(\bar x,\bar y)$ in direction $u$} whenever
			the following condition holds:
			for every sequences $\{(x_k,y_k)\}_{k\in\N}\subset\gph\Phi$,
			$\{x_k^*\}_{k\in\N}\subset\mathbb X$, and $\{\lambda_k\}_{k\in\N}\subset\mathbb Y$ 
			as well as $x^*\in\mathbb X$ and $y^*\in\mathbb Y$ satisfying 
			$x_k\notin\Phi^{-1}(\bar y)$, $y_k\neq\bar y$, and 
			$x_k^*\in \widehat{D}^*\Phi(x_k,y_k)(\lambda_k)$
			for each $k\in\N$ as well as the convergences
			\begin{equation}\label{eq:convergences_directional_asymptotic_regularity}
				\begin{aligned}
				x_k&\to\bar x,&
				\qquad
				y_k&\to\bar y,&
				\qquad
				x_k^*&\to x^*,&
				\\
				\frac{x_k-\bar x}{\norm{x_k-\bar x}}&\to u,&
				\qquad
				\frac{y_k-\bar y}{\norm{x_k-\bar x}}&\to 0,&
				\qquad  
				\norm{\lambda_k}&\to\infty,&
				\\
				\frac{y_k-\bar y}{\norm{y_k-\bar y}}-\frac{\lambda_k}{\norm{\lambda_k}}
				&\to 0,& \frac{\norm{y_k-\bar y}}{\norm{x_k-\bar x}}\lambda_k&\to y^*,& &&
				\end{aligned}
			\end{equation}
			we find $x^*\in\Im D^*\Phi(\bar x,\bar y)$.
		\item For the fixed direction $u\in\mathbb S_{\mathbb X}$, 
			$\Phi$ is said to be \emph{strongly asymptotically regular 
			at $(\bar x,\bar y)$ in direction $u$} whenever
			the following condition holds:
			for every sequences $\{(x_k,y_k)\}_{k\in\N}\subset\gph\Phi$,
			$\{x_k^*\}_{k\in\N}\subset\mathbb X$, and $\{\lambda_k\}_{k\in\N}\subset\mathbb Y$ 
			as well as $x^*\in\mathbb X$ and $y^*\in\mathbb Y$ satisfying 
			$x_k\notin\Phi^{-1}(\bar y)$, $y_k\neq\bar y$, 
			and $x_k^*\in \widehat{D}^*\Phi(x_k,y_k)(\lambda_k)$
			for each $k\in\N$ as well as the convergences
			\eqref{eq:convergences_directional_asymptotic_regularity}, 
			we have $x^*\in \Im D^*\Phi((\bar x,\bar y);(u,0))$.
	\end{enumerate}
\end{definition}

Let us briefly note that asymptotic regularity of a set-valued mapping 
$\Phi\colon\mathbb X\tto\mathbb Y$ at some point $(\bar x,0)\in\gph\Phi$ in the sense
of \cref{def:asymptotic_regularity} equals AM-regularity of the set $\Phi^{-1}(0)$
at $\bar x$ mentioned in \cite[Remark~3.17]{Mehlitz2020a}, see \cref{prop:asymptotic_regularity_via_limiting_tools} as well.
The concepts of directional asymptotic regularity from \cref{def:asymptotic_regularity} are new.

In the subsequent remark, we summarize some obvious relations between the different concepts
from \cref{def:asymptotic_regularity}.
\begin{remark}\label{rem:relations_asymptotic_regularity}
	Let $(\bar x,\bar y)\in\gph\Phi$ be fixed.
	Then the following assertions hold.
	\begin{enumerate}
		\item Let $u\in\mathbb S_{\mathbb X}$ be arbitrarily chosen.
			If $\Phi$ is strongly asymptotically regular at $(\bar x,\bar y)$ in direction
			$u$, it is asymptotically regular at $(\bar x,\bar y)$ in direction $u$.
		\item\label{item:dir_asymp_reg_equals_asymp_reg} 
			If $\Phi$ is asymptotically regular at $(\bar x,\bar y)$, then
			it is asymptotically regular at $(\bar x,\bar y)$ 
			in each direction from $\mathbb S_{\mathbb X}$.
	\end{enumerate}
\end{remark}

We note that strong asymptotic regularity in each unit direction is indeed not related to asymptotic regularity. 
On the one hand, the subsequently stated example, taken from \cite[Example~3.15]{Mehlitz2020a},
shows that asymptotic regularity does not imply strong asymptotic regularity in each unit direction.
On the other hand, \cref{ex:FOSCMS_but_not_asymptotically_regular}
illustrates that strong asymptotic regularity in each unit direction does not
yield asymptotic regularity.
\begin{example}\label{ex:asymptotic_regularity_but_not_strong}
	We consider $\Phi\colon\R\tto\R$ given by
	\[
		\forall x\in\R\colon\quad
		\Phi(x)
		:=
		\begin{cases}
			\R		&\text{if }x\leq 0,\\
			[x^2,\infty)	&\text{if }x>0
		\end{cases}
	\]
	at $(\bar x,\bar y):=(0,0)$. 
	It is demonstrated in \cite[Example~3.15]{Mehlitz2020a} that $\Phi$ is 
	asymptotically regular at $(\bar x,\bar y)$.
	We find $\mathcal T_{\gph\Phi}(\bar x,\bar y)=\{(u,v)\in\R^2\,|\,u\leq 0\,\lor\,v\geq 0\}$
	so $(\pm 1,0)\in\mathcal T_{\gph\Phi}(\bar x,\bar y)$.
	Let us consider $u:=1$.
	Then we find $\Im D^*\Phi((\bar x,\bar y);(u,0))=\{0\}$.
	Taking $x^*:=1$, $y^*:=1/2$, as well as
	\[
		\forall k\in\N\colon\quad
		x_k:=\frac1k,
		\qquad
		y_k:=\frac1{k^2},
		\qquad
		x_k^*:=1,
		\qquad
		\lambda_k:=\frac k2,
	\]
	we have $x_k^*\in\widehat D^*\Phi(x_k,y_k)(\lambda_k)$ for all $k\in\N$
	as well as the convergences \eqref{eq:convergences_directional_asymptotic_regularity}.
	However, due to $x_k^*\to x^*\notin \Im D^*\Phi((\bar x,\bar y);(u,0))$, $\Phi$ is
	not strongly asymptotically regular at $(\bar x,\bar y)$ in direction $u$.
\end{example}

Combining \cref{thm:directional_asymptotic_stationarity} with the concepts from \cref{def:asymptotic_regularity}, 
we immediately obtain the following result due to local boundedness of the
regular subdifferential of Lipschitzian functions, see e.g.\
\cite[Theorem~1.22]{Mordukhovich2018}.
\begin{corollary}\label{cor:M_stationarity_via_directional_asymptotic_regularity}
 	Let $\bar x\in\mathcal F$ be a local minimizer of \eqref{eq:nonsmooth_problem} 
 	such that, for each critical direction $u\in\mathbb S_{\mathbb X}$ 
 	for \eqref{eq:nonsmooth_problem} at $\bar x$, $\Phi$ is
	asymptotically regular at $(\bar x,\bar y)$ in direction $u$.
	Then $\bar x$ is M-stationary.
\end{corollary}

In the light of 
\cref{rem:relations_asymptotic_regularity}\,\ref{item:dir_asymp_reg_equals_asymp_reg}, 
our result from
\cref{cor:M_stationarity_via_directional_asymptotic_regularity} improves
\cite[Theorem~3.9]{Mehlitz2020a} by a directional refinement of the constraint qualification
since it suffices to check asymptotic regularity with respect to particular directions.

We point out that, unlike typical constraint qualifications, (directional) asymptotic regularity
allows the existence of sequences 
satisfying \eqref{eq:convergences_directional_asymptotic_regularity}
as long as the limit $x^*$ is included in $\Im D^*\Phi(\bar x,\bar y)$ 
which is enough for M-stationarity.

For the purpose of completeness, 
we show that the notions from \cref{def:asymptotic_regularity} can be stated
in terms of the limiting coderivative completely.

\begin{proposition}\label{prop:asymptotic_regularity_via_limiting_tools}
	\cref{def:asymptotic_regularity} can be equivalently formulated 
	with $x_k^*\in\widehat{D}^*\Phi(x_k,y_k)(\lambda_k)$ replaced by 
	$x_k^*\in D^*\Phi(x_k,y_k)(\lambda_k)$.
\end{proposition}
\begin{proof}
	For non-directional asymptotic regularity 
	the proof is standard and follows from a simple diagonal sequence argument.
	The proof for strong directional asymptotic regularity parallels the one 
	for directional asymptotic regularity which is presented below.

	Since one implication is clear by definition of the
			regular and limiting coderivative, we only show the other one.
			Thus, let us fix sequences $\{(x_k,y_k)\}_{k\in\N}\subset\gph\Phi$,
			$\{x_k^*\}_{k\in\N}\subset\mathbb X$, and 
			$\{\lambda_k\}_{k\in\N}\subset\mathbb Y$ as well as $x^*\in\mathbb X$ 
			and $y^*\in\mathbb Y$ satisfying 
			$x_k\notin\Phi^{-1}(\bar y)$, $y_k\neq\bar y$, 
			and $x_k^*\in D^*\Phi(x_k,y_k)(\lambda_k)$
			for each $k\in\N$ as well as the convergences
			\eqref{eq:convergences_directional_asymptotic_regularity}.
			For each $k\in\N$, we find sequences 
			$\{(x_{k,\ell},y_{k,\ell})\}_{\ell\in\N}\subset\gph\Phi$,
			$\{x_{k,\ell}^*\}_{\ell\in\N}\subset\mathbb X$, 
			and $\{\lambda_{k,\ell}\}_{\ell\in\N}\subset\mathbb Y$ with
			$x_{k,\ell}\to x_k$, $x_{k,\ell}^*\to x_k^*$, $y_{k,\ell}\to y_k$, 
			and $\lambda_{k,\ell}\to \lambda_k$
			as $\ell\to\infty$ as well as 
			$x_{k,\ell}^*\in \widehat D^*\Phi(x_{k,\ell},y_{k,\ell})(\lambda_{k,\ell})$
			for each $\ell\in\N$. 
			Observing that $\Phi^{-1}(\bar y)$ is closed, 
			its complement is open so that $x_{k,\ell}\notin\Phi^{-1}(\bar y)$
			holds for sufficiently large $\ell\in\N$. 
			Furthermore, since $\norm{x_k-\bar x}>0$ and $\norm{y_k-\bar y}>0$ are valid, 
			we can choose
			an index $\ell(k)\in\N$ so large such that the estimates
			\[
				\begin{aligned}
				\nnorm{x_{k,\ell(k)}-x_k}&<\frac{1}{k}\norm{x_k-\bar x},&\quad
				\nnorm{x_{k,\ell(k)}^*-x_k^*}&<\frac{1}{k},&
				\\
				\nnorm{y_{k,\ell(k)}-y_k}&
					<\frac{1}{k}\norm{y_k-\bar y},&\quad
				\nnorm{\lambda_{k,\ell(k)}-\lambda_k}&<\frac1k&
				\end{aligned}
			\]
			and $x_{k,\ell(k)}\notin\Phi^{-1}(\bar y)$ as well as $y_{k,\ell(k)}\neq\bar y$ 
			are valid. 
			For each $k\in\N$, we set $\tilde x_k:=x_{k,\ell(k)}$, 
			$\tilde x_k^*:=x_{k,\ell(k)}^*$,
			$\tilde y_k:=y_{k,\ell(k)}$, and $\tilde\lambda_k:=\lambda_{k,\ell(k)}$.
			Clearly, we have $\tilde x_k\to\bar x$, $\tilde y_k\to\bar y$, 
			$\tilde x_k^*\to x^*$, $\nnorm{\tilde \lambda_k}\to\infty$,
			$\{(\tilde x_k,\tilde y_k)\}_{k\in\N}\subset\gph\Phi$, 
			and $\tilde x_k\notin\Phi^{-1}(\bar y)$, $\tilde y_k\neq\bar y$,
			as well as $\tilde x_k^*\in\widehat D^*\Phi(\tilde x_k,\tilde y_k)(\tilde\lambda_k)$
			for each $k\in\N$ by construction.	
			Furthermore, we find
			\[
				\norm{\tilde x_k-\bar x}
				\geq
				\norm{x_k-\bar x}-\norm{\tilde x_k-x_k}
				\geq
				\frac{k-1}{k}\norm{x_k-\bar x}
			\]
			for each $k\in\N$. 
			With the above estimates at hand, we obtain
			\begin{align*}
				\norm{
					\frac{x_k-\bar x}{\nnorm{x_k-\bar x}}
					-
					\frac{\tilde x_k-\bar x}{\norm{\tilde x_k-\bar x}}
				}
				&=
				\norm{
					\frac{x_k-\tilde x_k}{\norm{x_k-\bar x}}
					+
					(\tilde x_k-\bar x)
					\left(
						\frac{1}{\nnorm{x_k-\bar x}}-\frac{1}{\nnorm{\tilde x_k-\bar x}}
					\right)
				}
				\\
				&
				\leq 
				\frac{\norm{x_k-\tilde x_k}}{\norm{x_k-\bar x}}
				+
				\frac{\norm{\tilde x_k-\bar x}\norm{x_k-\tilde x_k}}
					{\norm{x_k-\bar x}\norm{\tilde x_k-\bar x}}		
				\leq 
				\frac{2}{k}		
			\end{align*}
			and
			\begin{equation}\label{eq:estimate_on_norm_quotients}
				\begin{aligned}
				\norm{
					\frac{y_k-\bar y}{\nnorm{x_k-\bar x}}
					-
					\frac{\tilde y_k-\bar y}{\norm{\tilde x_k-\bar x}}
				}
				&=
				\norm{
					\frac{y_k-\tilde y_k}{\norm{x_k-\bar x}}
					+
					(\tilde y_k-\bar y)
					\left(
						\frac{1}{\nnorm{x_k-\bar x}}-\frac{1}{\nnorm{\tilde x_k-\bar x}}
					\right)
				}
				\\
				&
				\leq
				\frac{\norm{y_k-\tilde y_k}}{\norm{x_k-\bar x}}
				+
				\frac{\norm{\tilde y_k-\bar y}{\norm{x_k-\tilde x_k}}}
					{\norm{x_k-\bar x}\norm{\tilde x_k-\bar x}}
				\\
				&\leq
				\frac1k\frac{\norm{y_k-\bar y}}{\norm{x_k-\bar x}}
				+
				\frac{1}{k-1}\frac{\norm{\tilde y_k-y_k}+\norm{y_k-\bar y}}{\norm{x_k-\bar x}}
				\\
				&\leq
				\left(\frac{1}{k} + \frac{1}{k(k-1)} + \frac{1}{k-1}\right)
					\frac{\norm{y_k-\bar y}}{\norm{x_k-\bar x}}
				\\
				&
				=
				\frac{2}{k-1}\frac{\norm{y_k-\bar y}}{\norm{x_k-\bar x}},
				\end{aligned}
			\end{equation}
			so that, with the aid of \eqref{eq:convergences_directional_asymptotic_regularity},
			we find
			$(\tilde x_k-\bar x)/\norm{\tilde x_k-\bar x}\to u$ 
			and $(\tilde y_k-\bar y)/\norm{\tilde x_k-\bar x}\to 0$.
			With the aid of \eqref{eq:estimate_on_norm_quotients},
			\begin{align*}
				\norm{
					\frac{\nnorm{\tilde y_k-\bar y}}{\nnorm{\tilde x_k-\bar x}}
					\tilde\lambda_k
					-
					\frac{\norm{y_k-\bar y}}{\norm{x_k-\bar x}}\lambda_k
				}
				&
				\leq
				\frac{\nnorm{\tilde y_k-\bar y}}{\nnorm{\tilde x_k-\bar x}}
					\nnorm{\tilde\lambda_k-\lambda_k}
				+
				\left|
					\frac{\nnorm{\tilde y_k-\bar y}}{\nnorm{\tilde x_k-\bar x}}
					-
					\frac{\norm{y_k-\bar y}}{\norm{x_k-\bar x}}
				\right|
					\norm{\lambda_k}
				\\
				&
				\leq
				\frac1k\frac{\nnorm{\tilde y_k-\bar y}}{\nnorm{\tilde x_k-\bar x}}
				+
				\frac{2}{k-1}\frac{\norm{y_k-\bar y}}{\norm{x_k-\bar y}}\norm{\lambda_k}
			\end{align*}
			is obtained, which gives 
			$\tilde\lambda_k\nnorm{\tilde y_k-\bar y}/\nnorm{\tilde x_k-\bar x}\to y^*$.
			Similar as above, we find
			\begin{align*}
				\norm{\frac{\tilde y_k-\bar y}{\norm{\tilde y_k-\bar y}}
					-\frac{y_k-\bar y}{\norm{y_k-\bar y}}}
				\leq\frac{2}{k}
			\end{align*}
			and
			\begin{align*}
				\nnorm{
					\tilde\lambda_k/\nnorm{\tilde\lambda_k}
					-
					\lambda_k/\norm{\lambda_k}
					}
				\leq
				2\nnorm{\lambda_k-\tilde\lambda_k}/\norm{\lambda_k}
				\leq
				2/(k\norm{\lambda_k}),
			\end{align*}
			so that \eqref{eq:convergences_directional_asymptotic_regularity} gives us
			\begin{align*}
				\lim\limits_{k\to\infty}
				\left(\frac{\tilde y_k-\bar y}{\norm{\tilde y_k-\bar y}}
					-\frac{\tilde\lambda_k}{\nnorm{\tilde\lambda_k}}\right)
				=
				\lim\limits_{k\to\infty}
				\left(\frac{y_k-\bar y}{\norm{y_k-\bar y}}
					-\frac{\lambda_k}{\norm{\lambda_k}}\right)
				=
				0.
			\end{align*}
			Now, since $\Phi$ is asymptotically regular at $(\bar x,\bar y)$ in direction $u$, 
			we obtain $x^*\in\Im D^*\Phi(\bar x,\bar y)$.
\end{proof}

Since (directional) asymptotic regularity (w.r.t.\ all critical unit directions) yields
M-stationarity of a local minimizer 
by \cref{cor:M_stationarity_via_directional_asymptotic_regularity},
in the remaining part of the paper, we put it into context of other
common assumptions that work as a constraint qualification
for M-stationarity associated with problem \eqref{eq:nonsmooth_problem}.
Let us clarify here some rather simple or known connections.
\begin{enumerate}
 \item\label{item:polyhedrality} 
 	A polyhedral mapping is asymptotically regular at each point of its graph.
 \item\label{item:metric_regularity} 
 	Metric regularity implies asymptotic regularity.
 \item\label{item:strong_metric_subregularity} 
 	Strong metric subregularity implies asymptotic regularity.
 \item\label{item:FOSCMS_does_not_imply_asymptotic_regularity} 
 	FOSCMS does not imply asymptotic regularity,
 	but it implies strong asymptotic regularity in each unit direction.
 \item\label{item:metric_subregularity_not_sufficient_for_asymptotic_regularity}
 	Metric subregularity does not imply asymptotic regularity in each unit direction.
 	However, if the map of interest is metrically subregular 
 	at every point of its graph near the reference point with a \emph{uniform} constant,
 	then strong asymptotic regularity in each unit direction follows.
 \item\label{item:asymptotic_regularity_does_not_yield_exact_penalty}
 	Neither asymptotic regularity nor strong directional asymptotic regularity yields
 	the directional exact penalty property of \cref{lem:directional_exact_penalization}.
\end{enumerate}

Statements~\ref{item:polyhedrality} and~\ref{item:metric_regularity} 
were shown in \cite[Theorems 3.10 and 3.12]{Mehlitz2020a}.
Let us now argue that strong metric subregularity (the ``inverse'' property associated with
isolated calmness), see \cite{DontchevRockafellar2014}, 
also implies asymptotic regularity at the point.
This follows easily from the discussion above
\cite[Corollary~4.6]{BenkoMehlitz2020}, which yields
that the domain of the limiting coderivative, at the point 
where the mapping is isolatedly calm, is the whole space.
Equivalently, the range of the limiting coderivative, at the point 
where the mapping is strongly metrically subregular, is the whole space 
and asymptotic regularity thus follows trivially.
Thus, statement~\ref{item:strong_metric_subregularity} follows.

Regarding~\ref{item:FOSCMS_does_not_imply_asymptotic_regularity}, the fact that FOSCMS implies
strong asymptotic regularity in each unit direction easily follows
by similar arguments that show that
metric regularity implies asymptotic regularity, 
see \cite[Lemma~3.11, Theorem~3.12]{Mehlitz2020a}.
Indeed, let us fix $(\bar x,\bar y)\in\gph\Phi$ and $u\in\mathbb S_{\mathbb X}$ such that 
$\ker D^*\Phi((\bar x,\bar y);(u,0))=\{0\}$. 
Furthermore, choose sequences $\{(x_k,y_k)\}_{k\in\N}\subset\gph\Phi$,
$\{x_k^*\}_{k\in\N}\subset\mathbb X$, and $\{\lambda_k\}_{k\in\N}\subset\mathbb Y$ 
as well as $x^*\in\mathbb X$ and $y^*\in\mathbb Y$
satisfying $x_k\notin\Phi^{-1}(\bar y)$, $y_k\neq\bar y$, and 
$x_k^*\in \widehat{D}^*\Phi(x_k,y_k)(\lambda_k)$ for all
$k\in\N$ as well as the convergences \eqref{eq:convergences_directional_asymptotic_regularity}.
Then we also have 
$x_k^*/\norm{\lambda_k}\in \widehat{D}^*\Phi(x_k,y_k)(\lambda_k/\norm{\lambda_k})$ 
for each $k\in\N$,
and taking the limit $k\to\infty$ along a suitably chosen subsequence, 
we end up with $0\in D^*\Phi((\bar x,\bar y);(u,0))(\lambda)$ 
for some $\lambda\in\mathbb S_{\mathbb Y}$ which is a contradiction.
Hence, such sequences cannot exist and $\Phi$ is strongly asymptotically 
regular at $(\bar x,\bar y)$ in direction $u$.

The following example shows that FOSCMS
does not imply asymptotic regularity.
\begin{example}\label{ex:FOSCMS_but_not_asymptotically_regular}
    Let $\Phi\colon \R \tto \R$ be given by
    \begin{equation*}
        \forall x\in\R\colon\quad
        \Phi(x):= 
        \begin{cases}
             [x,\infty) & \text{if } x \leq 0, \\
             \left[\frac{1}{k} - \frac{1}{k}\left(x - \frac{1}{k}\right),\infty\right)
             	& \text{if } x \in \left(\frac{1}{k+1},\frac{1}{k}\right] \text{ for some }k\in\N,\\
             \varnothing	&\text{otherwise.}
        \end{cases}
    \end{equation*}
    Then $\{(1/k,1/k)\}_{k\in\N}\subset\gph \Phi$ converges to
    $(\bar x,\bar y):=(0,0)$ and
    \begin{equation*}
        \mathcal N_{\gph \Phi}(1/k,1/k) =
        \{(x^*,y^*)\in\R^2 \,\vert\, y^* \leq 0, y^* \leq k x^*\}
    \end{equation*}
    is valid showing that $\Im D^*\Phi(1/k,1/k) = \R$ is valid for all $k\in\N$.
    On the other hand, we have
    \begin{equation*}
        \mathcal N_{\gph \Phi}(0,0) =
        \{(x^*,y^*)\in\R^2 \,\vert\, x^* \geq 0, y^* \leq 0\},
    \end{equation*}
    and thus $\Im D^*\Phi(0,0) = \R_+$.
    This means that $\Phi$ is not asymptotically regular at $(\bar x,\bar y)$.
    
    On the other hand, we find
    \[
    	\mathcal T_{\gph\Phi}(\bar x,\bar y)=\{(u,v)\in\R^2\,|\,u\leq v\}.
    \]
    Each pair $(u,0)\in \mathcal T_{\gph\Phi}(\bar x,\bar y)$ with $u\neq 0$ satisfies
    $u<0$, i.e., the direction $(u,0)$ points into the interior of $\gph\Phi$. 
    Thus, we have $\mathcal N_{\gph\Phi}((\bar x,\bar y),(u,0))=\{(0,0)\}$
    which shows that FOSCMS is valid.
\end{example}

Regarding~\ref{item:metric_subregularity_not_sufficient_for_asymptotic_regularity}, 
let us fix $(\bar x,\bar y)\in\gph\Phi$ and note that metric subregularity of $\Phi$ 
on a neighborhood of $(\bar x,\bar y)$ (restricted to $\gph\Phi$) with a uniform constant
$\kappa>0$ is clearly milder
than metric regularity at $(\bar x,\bar y)$ 
since it is automatically satisfied e.g.\ by polyhedral mappings.
To see that it implies asymptotic regularity, consider
sequences $\{(x_k,y_k)\}_{k\in\N}\subset\gph\Phi$, $\{x_k^*\}_{k\in\N}\subset\mathbb X$,
and $\{\lambda_k\}_{k\in\N}\subset\mathbb Y$ as well as $x^*\in\mathbb X$ 
and $y^*\in\mathbb Y$ satisfying 
$x_k^*\in \widehat{D}^*\Phi(x_k,y_k)(\lambda_k)$ for each $k\in\N$ 
and the convergences \eqref{eq:convergences_directional_asymptotic_regularity}
for some unit direction $u\in\mathbb S_{\mathbb X}$.
Due to \cite[Theorem 3.2]{BenkoMehlitz2020} and 
$-x_k^*\in \dom\widehat{D}^*\Phi^{-1}(y_k,x_k)$, we find
$x_k^*\in \widehat{\mathcal N}_{\Phi^{-1}(y_k)}(x_k) \subset \mathcal N_{\Phi^{-1}(y_k)}(x_k)$
for each $k\in\N$.
Furthermore, \cite[Theorem 3.2]{BenkoMehlitz2020} also gives the existence of 
$\tilde\lambda_k\in\mathbb Y$
with $\nnorm{\tilde\lambda_k} \leq \kappa \norm{x_k^*}$ 
and $x_k^*\in D^*\Phi(x_k,y_k)(\tilde\lambda_k)$.
Noting that $\{x_k^*\}_{k\in\N}$ converges, this shows that there is 
a limit point $\lambda\in\mathbb Y$ 
of $\{\tilde\lambda_k\}_{k\in\N}$ which satisfies 
$x^*\in D^*\Phi((\bar x,\bar y);(u,0))(\lambda)$
by robustness of the directional limiting coderivative 
which can be distilled from \cref{lem:robustness_directional_limiting_normals}.
Hence, $\Phi$ is strongly asymptotically regular at $(\bar x,\bar y)$ in direction $u$.
Note that for the above arguments to work, we only need uniform metric subregularity 
along all sequences $\{(x_k,y_k)\}_{k\in\N}\subset\gph\Phi$
converging to $(\bar x, \bar y)$ from direction $(u,0)$.

On the other hand, the following example shows that
metric subregularity in the neighborhood of the point of interest
does not imply asymptotic regularity in each unit direction.

\begin{example}\label{ex:metric_subregularity_vs_asymptotic_regularity}
	We consider the mapping $\Phi\colon\R\tto\R$ given by
	\[
		\forall x\in\R\colon\quad
		\Phi(x):=\{0,x^2\}.
	\]
	Due to $\Phi^{-1}(0)=\R$, $\Phi$ is metrically subregular at all points $(x,0)$ where $x\in\R$
	is arbitrary.
	Furthermore, at all points $(x,x^2)$ where $x\neq 0$ holds, the Mordukhovich criterion
	shows that $\Phi$ is metrically regular. 
	Thus, $\Phi$ is metrically subregular at each point of its graph.
	Note that the moduli of metric subregularity tend to $\infty$ along the points 
	$(t,t^2)$ and $(-t,t^2)$ as $t\searrow 0$.
	
	Let us consider the point $(\bar x,\bar y):=(0,0)$ where we have 
	$\mathcal N_{\gph\Phi}(\bar x,\bar y)=\{0\}\times\R$
	and, thus, $\Im D^*\Phi(\bar x,\bar y)=\{0\}$. 
	Choosing $x^*:=1$, $y^*:=1/2$, as well as
	\[
		\forall k\in\N\colon\quad
		x_k:=\frac1k,\qquad y_k:=\frac{1}{k^2},\qquad x_k^*:=1,\qquad \lambda_k:=\frac{k}{2},
	\]
	we have $x_k^*\in\widehat D^*\Phi(x_k,y_k)(\lambda_k)$ for all $k\in\N$ as well as
	the convergences \eqref{eq:convergences_directional_asymptotic_regularity} for $u:=1$.
	Due to $x_k^*\to x^*\notin\Im D^*\Phi(\bar x,\bar y)$, $\Phi$ is not asymptotically regular
	at $(\bar x,\bar y)$ in direction $u$.
\end{example}

Finally, let us address item~\ref{item:asymptotic_regularity_does_not_yield_exact_penalty}
with the aid of an example.
\begin{example}\label{ex:asymptotic_regularity_and_no_exact_penalty}
	Let us define $\varphi\colon\R\to\R$ and $\Phi\colon\R\tto\R$ by means of
	\[
		\forall x\in\R\colon\quad
		\varphi(x):=-x,
		\qquad
		\Phi(x)
		:=
		\begin{cases}
			\R	&\text{if }x\leq 0,\\
			[x^2,\infty)	&\text{if }x=\frac1k\text{ for some }k\in\N,\\
			\varnothing	&\text{otherwise.}
		\end{cases}
	\]
	Furthermore, we fix $\bar y:=0$.
	One can easily check that $\bar x:=0$ is the uniquely determined global
	minimizer of the associated problem \eqref{eq:nonsmooth_problem}.
	Furthermore, we have 
	$\Im D^*\Phi(\bar x,\bar y)=\Im D^*\Phi((\bar x,\bar y);(1,0))=\R$ which
	shows that $\Phi$ is asymptotically regular at $(\bar x,\bar y)$ as well
	as strongly asymptotically regular at $(\bar x,\bar y)$ in direction $1$.
	Furthermore, it is obvious that $\Phi$ is strongly asymptotically regular
	at $(\bar x,\bar y)$ in direction $-1$.
	Finally, let us mention that $\Phi$ fails to be metrically subregular at
	$(\bar x,\bar y)$ in direction $1$.
	
	Now, define $x_k:=1/k$ for each $k\in\N$ and observe that for each constant
	$C>0$ and sufficiently large $k\in\N$, we have 
	$\varphi(x_k)+C\,\dist(\bar y,\Phi(x_k))=-1/k+C/k^2<0=\varphi(\bar x)$, i.e., 
	$\bar x$ is not a local minimizer of \eqref{eq:directionally_penalized_problem}
	for any choice of $C>0$, $\varepsilon>0$, $\delta>0$, and $u:=1$.	
\end{example}

\section{Directional pseudo- and quasi-normality}\label{sec:pseudo_quasi_normality}

In this section, we connect asymptotic regularity with the notions of pseudo- and quasi-normality.
Note that the latter concepts have been introduced for
standard nonlinear programs in 
\cite{BertsekasOzdaglar2002,Hestenes1975},
and reasonable generalizations to more general
geometric constraints have been established in \cite{GuoYeZhang2013}.
Furthermore, problem-tailored notions of these conditions have
been coined e.g.\ for so-called cardinality-, complementarity-, and switching-constrained
optimization problems, see \cite{KanzowRaharjaSchwartz2021b,KanzowSchwartz2010,LiangYe2021}.
Let us point out that these conditions are comparatively
mild constraint qualifications and sufficient for the
presence of metric subregularity of the underlying feasibility mapping 
which equals the so-called error bound property,
see e.g.\ \cite[Theorem~5.2]{GuoYeZhang2013}.
Here, we extend pseudo- and quasi-normality
from the common setting of geometric constraint
systems to arbitrary set-valued mappings and comment
on the qualitative properties of these conditions.
Naturally, we aim for directional versions of these concepts, which,
in the setting of geometric constraints,
were recently introduced in \cite{BaiYeZhang2019} and further explored
in \cite{BenkoCervinkaHoheisel2019}.
Furthermore, we briefly discuss directional pseudo- and quasi-normality
in the context of equilibrium-constrained optimization.

\subsection{Pseudo- and quasi-normality for set-valued mappings}

The definition below introduces the notions of our interest.

\begin{definition}\label{def:direction_quasi_pseudo_normality}
	Fix $(\bar x,\bar y)\in\gph\Phi$ and a direction $u\in\mathbb S_{\mathbb X}$.
	\begin{enumerate}
		\item We say that \emph{pseudo-normality in direction $u$} holds at $(\bar x,\bar y)$
			if there does not exist a nonzero vector 
			$\lambda\in\ker D^*\Phi((\bar x,\bar y);(u,0))$
			satisfying the following condition:
			there are sequences $\{(x_k,y_k)\}_{k\in\N}\subset\gph\Phi$ 
			with $x_k\neq\bar x$ for all $k\in\N$ and
			$\{\lambda_k\}_{k\in\N}\subset\mathbb Y$, $\{\eta_k\}_{k\in\N}\subset\mathbb X$,
			such that
			\begin{equation}\label{eq:convergences_definition_quasi_normality}
				\begin{aligned}
				x_k&\to\bar x,&
				\qquad
				y_k&\to\bar y,&
				\qquad
				\lambda_k&\to\lambda,&
				\\
				\eta_k&\to 0,&
				\qquad
				\frac{x_k-\bar x}{\norm{x_k-\bar x}}&\to u,&
				\qquad
				\frac{y_k-\bar y}{\norm{x_k-\bar x}}&\to 0,&
				\end{aligned}
			\end{equation}
			and $\eta_k\in \widehat{D}^*\Phi(x_k,y_k)(\lambda_k)$ 
			as well as $\dual{\lambda}{y_k-\bar y}>0$
			for all $k\in\N$.
		\item Let $\mathcal E:=\{e_1,\ldots,e_m\}\subset\mathbb Y$ be an orthonormal basis of
			$\mathbb Y$. We say that \emph{quasi-normality in direction $u$} holds at $(\bar x,\bar y)$
			w.r.t.\ $\mathcal E$ if there does not exist a nonzero vector 
			$\lambda\in\ker D^*\Phi((\bar x,\bar y);(u,0))$
			satisfying the following condition:
			there are sequences $\{(x_k,y_k)\}_{k\in\N}\subset\gph\Phi$ 
			with $x_k\neq\bar x$ for all $k\in\N$ and
			$\{\lambda_k\}_{k\in\N}\subset\mathbb Y$, $\{\eta_k\}_{k\in\N}\subset\mathbb X$,
			such that we have the convergences from
			\eqref{eq:convergences_definition_quasi_normality}
			and, for all $k\in\N$ and $i\in\{1,\ldots,m\}$, 
			$\eta_k\in \widehat{D}^*\Phi(x_k,y_k)(\lambda_k)$ 
			as well as $\dual{\lambda}{e_i}\dual{y_k-\bar y}{e_i}>0$
			if $\dual{\lambda}{e_i}\neq 0$.
	\end{enumerate}
\end{definition}

In case where the canonical basis is chosen in $\mathbb Y:=\R^m$, the
above concept of quasi-normality is a direct generalization of the original
notion from \cite{BertsekasOzdaglar2002} which was coined for standard
nonlinear problems and neglected directional information.
Let us just mention that a reasonable, basis-independent definition of quasi-normality 
would require that there exists some basis w.r.t.\ which the mapping of interest is quasi-normal,
see also \cref{thm:quasi_normality_yields_asymptotic_regularity}.

Note that the sequence $\{y_k\}_{k\in\N}$ in the definition of directional pseudo- and
quasi-normality needs to satisfy $y_k\neq\bar y$ for all $k\in\N$. In the definition
of directional pseudo-normality, this is clear from $\dual{\lambda}{y_k-\bar y}>0$ for all
$k\in\N$. Furthermore, in the definition of directional quasi-normality, observe that
$\lambda\neq 0$ implies the existence of $j\in\{1,\ldots,m\}$ such that 
$\dual{\lambda}{e_j}\neq 0$ holds, so that $\dual{y_k-\bar y}{e_j}\neq 0$ is necessary
for each $k\in\N$.

In the following lemma, we show the precise relation between directional pseudo- and quasi-normality.

\begin{lemma}\label{lem:pseudo_vs_quasi_normality}
	Fix $(\bar x,\bar y)\in\gph\Phi$ and some direction $u\in\mathbb S_{\mathbb X}$.
	Then $\Phi$ is pseudo-normal at $(\bar x,\bar y)$ in direction $u$ if and only
	if $\Phi$ is quasi-normal at $(\bar x,\bar y)$ in direction $u$ w.r.t.\
	each orthonormal basis of $\mathbb Y$.
\end{lemma}
\begin{proof}
	$[\Longrightarrow]$ Let $\Phi$ be pseudo-normal at $(\bar x,\bar y)$ in direction $u$,
		let $\mathcal E:=\{e_1,\ldots,e_m\}\subset\mathbb Y$ be an orthonormal basis of
		$\mathbb Y$, and pick $\lambda\in\ker D^*\Phi((\bar x,\bar y);(u,0))$
		as well as sequences $\{(x_k,y_k)\}_{k\in\N}\subset\gph\Phi$ 
		with $x_k\neq\bar x$ for all $k\in\N$ and
		$\{\lambda_k\}_{k\in\N}\subset\mathbb Y$, $\{\eta_k\}_{k\in\N}\subset\mathbb X$,
		satisfying the convergences \eqref{eq:convergences_definition_quasi_normality}
		and, for all $k\in\N$ and $i\in\{1,\ldots,m\}$, 
		$\eta_k\in \widehat{D}^*\Phi(x_k,y_k)(\lambda_k)$ 
		as well as $\dual{\lambda}{e_i}\dual{y_k-\bar y}{e_i}>0$
		if $\dual{\lambda}{e_i}\neq 0$.
		Observing that we have
		\begin{align*}
			\dual{\lambda}{y_k-\bar y}
			&=
			\dual{\sum\nolimits_{i=1}^m\dual{\lambda}{e_i}e_i}{\sum\nolimits_{j=1}^m\dual{y_k-\bar y}{e_j}e_j}
			\\
			&=
			\sum\nolimits_{i=1}^m\sum\nolimits_{j=1}^m\dual{\lambda}{e_i}\dual{y_k-\bar y}{e_j}\dual{e_i}{e_j}
			\\
			&=
			\sum\nolimits_{i=1}^m\dual{\lambda}{e_i}\dual{y_k-\bar y}{e_i},
		\end{align*}
		validity of pseudo-normality at $(\bar x,\bar y)$ in direction $u$ gives $\lambda=0$,
		i.e., $\Phi$ is quasi-normal at $(\bar x,\bar y)$ in direction $u$ w.r.t.\ $\mathcal E$.\\
	$[\Longleftarrow]$ Assume that $\Phi$ is quasi-normal at $(\bar x,\bar y)$ in direction $u$ w.r.t.\
		each orthonormal basis of $\mathbb Y$.
		Suppose that $\Phi$ is not pseudo-normal at $(\bar x,\bar y)$ in direction $u$.
		Then we find some nonzero $\lambda\in\ker D^*\Phi((\bar x,\bar y);(u,0))$
		as well as sequences $\{(x_k,y_k)\}_{k\in\N}\subset\gph\Phi$ 
		with $x_k\neq\bar x$ for all $k\in\N$ and
		$\{\lambda_k\}_{k\in\N}\subset\mathbb Y$, $\{\eta_k\}_{k\in\N}\subset\mathbb X$,
		satisfying the convergences \eqref{eq:convergences_definition_quasi_normality}
		and	$\eta_k\in \widehat{D}^*\Phi(x_k,y_k)(\lambda_k)$ 
		as well as $\dual{\lambda}{y_k-\bar y}>0$ for all $k\in\N$.
		Noting that $\lambda$ does not vanish, we can construct an orthonormal basis
		$\mathcal E_\lambda:=\{e_1^\lambda,\ldots,e_m^\lambda\}$ of $\mathbb Y$ with
		$e_1^\lambda:=\lambda/\norm{\lambda}$. 
		Note that, for $i\in\{1,\ldots,m\}$, 
		we have $\ninnerprod{\lambda}{e_i^\lambda}\neq 0$ if and only
		if $i=1$ by construction of $\mathcal E_\lambda$. Furthermore, we find
		\begin{align*}
			\ninnerprod{\lambda}{e_1^\lambda}\ninnerprod{y_k-\bar y}{e_1^\lambda}
			=
			\norm{\lambda}\ninnerprod{\lambda/\norm{\lambda}}{y_k-\bar y}
			=
			\ninnerprod{\lambda}{y_k-\bar y}
			>
			0.
		\end{align*}
		This, however, contradicts quasi-normality of $\Phi$ at $(\bar x,\bar y)$ in direction $u$
		w.r.t.\ $\mathcal E_{\lambda}$.
\end{proof}

Let us note that \cite[Example~1]{BertsekasOzdaglar2002} shows in the non-directional
situation of standard nonlinear programming that pseudo-normality might be more 
restrictive than quasi-normality w.r.t.\ the canonical basis in $\R^m$.
On the other hand, due to \cref{lem:pseudo_vs_quasi_normality}, there must exist another 
basis such that quasi-normality w.r.t.\ this basis fails since pseudo-normality fails.
This depicts that validity of quasi-normality indeed may depend on the chosen basis.
In \cite{BaiYeZhang2019}, the authors define directional quasi-normality 
for geometric constraints in Euclidean spaces in componentwise fashion although this 
is somehow unclear in situations where the image space is different from $\R^m$. 
Exemplary, in the $\tfrac12m(m+1)$-dimensional space $\mathcal S_m$
of all real symmetric $m\times m$-matrices, the canonical basis, 
which seems to be associated with a componentwise calculus, comprises 
$\tfrac12(m-1)m$ matrices with precisely two nonzero entries. 
Our definition of quasi-normality from \cref{def:direction_quasi_pseudo_normality} 
gives some more freedom since the choice of the underlying basis allows to \emph{rotate}
the coordinate system.

Following the arguments in \cite[Section~3.2]{BenkoCervinkaHoheisel2019},
it also might be reasonable to define intermediate conditions bridging
pseudo- and quasi-normality. In the light of this paper, however, the concepts
from \cref{def:direction_quasi_pseudo_normality} are sufficient for
our purposes.

As the following theorem shows, directional quasi- and, thus, pseudo-normality
also serve as sufficient conditions for strong directional asymptotic regularity
and directional metric subregularity which explains our interest in these conditions.
Both statements follow once we clarify that pseudo- and quasi-normality
are in fact specifications of the multiplier sequential information in  \eqref{eq:convergences_directional_asymptotic_regularity},
namely $(y_k-\bar y)/\nnorm{y_k-\bar y} - \lambda_k/\nnorm{\lambda_k} \to 0$.

\begin{theorem}\label{thm:quasi_normality_yields_asymptotic_regularity}
  If $\Phi\colon\mathbb X \tto \mathbb Y$ is quasi-normal in direction
  $u \in \mathbb S_{\mathbb X}$ at $(\bar x, \bar y) \in \gph \Phi$
  w.r.t.\ some orthonormal basis $\mathcal E:=\{e_1,\ldots,e_m\}\subset\mathbb Y$ 
  of $\mathbb Y$, then it is also
  strongly asymptotically regular as well as metrically subregular in direction $u$ at $(\bar x, \bar y)$.
\end{theorem}
\begin{proof}
	Fix arbitrary sequences $\{(x_k,y_k)\}_{k\in\N}\subset\gph\Phi$,
	$\{x_k^*\}_{k\in\N}\subset\mathbb X$, and $\{\lambda_k\}_{k\in\N}\subset\mathbb Y$ 
	as well as $x^*\in\mathbb X$ and $y^*\in\mathbb Y$ satisfying 
	$x_k\notin\Phi^{-1}(\bar y)$, $y_k\neq\bar y$, and 
	$x_k^*\in \widehat{D}^*\Phi(x_k,y_k)(\lambda_k)$
	for each $k\in\N$ as well as the convergences
	\eqref{eq:convergences_directional_asymptotic_regularity}.
	Let us define $w_k:=(y_k-\bar y)/\norm{y_k-\bar y}$ 
	and $\tilde\lambda_k:=\lambda_k/\norm{\lambda_k}$ for each $k\in\N$.
	The requirements from \eqref{eq:convergences_directional_asymptotic_regularity} 
	imply that
	$\{w_k\}_{k\in\N}$ and  $\{\tilde\lambda_k\}_{k\in\N}$ converge, 
	along a subsequence (without relabeling),
	to the same nonvanishing limit which we will call $\lambda\in\mathbb S_{\mathbb Y}$.
	Moreover, given $i\in\{1,\ldots,m\}$ with $\dual{\lambda}{e_i} \neq 0$,
	for sufficiently large $k\in\N$, we get $\dual{w_{k}}{e_i} \neq 0 $ and
	\[
		0 
		< 
		\dual{\lambda}{e_i}\dual{ w_{k}}{e_i} 
		= 
		\dual{\lambda}{e_i}\dual{y_k - \bar y}{e_i} /\norm{y_k - \bar y}.
	\]
	Observing that we have $x_k^*/\norm{\lambda_k}\to 0$ from
	\eqref{eq:convergences_directional_asymptotic_regularity}, we find
	$\lambda\in\ker D^*\Phi((\bar x,\bar y);(u,0))$ by definition of the 
	directional limiting coderivative.
	This contradicts validity of quasi-normality of $\Phi$ at $(\bar x,\bar y)$
	in direction $u$ w.r.t.\ $\mathcal E$.
	Particularly, such sequences $\{(x_k,y_k)\}_{k\in\N}$, $\{x_k^*\}_{k\in\N}$, 
	and $\{\lambda_k\}_{k\in\N}$ cannot exist which means
	that $\Phi$ is strongly asymptotically regular in direction $u$ at $(\bar x,\bar y)$.

	The claim about metric subregularity now follows from \cite[Corollary 1]{Gfrerer2014a},
	since the only difference from quasi-normality is the requirement
	\[
	\dual{\lambda_k/\nnorm{\lambda_k}}{(y_k-\bar y)/\nnorm{y_k-\bar y}} \to 1
	\]
	which is the same as $(y_k-\bar y)/\nnorm{y_k-\bar y} - \lambda_k/\nnorm{\lambda_k} \to 0$
	as mentioned in the comments at the end of \cref{sec:model_program}.
\end{proof}

Relying on this result, \cref{lem:directional_M_stationarity_via_metric_subregularity}
yields that directional pseudo- and quasi-normality provide constraint qualifications
for \eqref{eq:nonsmooth_problem} which ensure validity of directional M-stationarity
at local minimizers.

We would like to point the reader's attention to the fact that non-directional versions of
pseudo- and quasi-normality are not comparable with the non-directional version of
asymptotic regularity. This has been observed in the context of standard nonlinear
programming, see \cite[Sections~4.3, 4.4]{AndreaniMartinezRamosSilva2016}.
The reason is that the standard version of asymptotic regularity makes no use of the multiplier information \eqref{eq:convergences_gamma=1_multiplier}.

Let us now also justify the terminology by showing that the new notions from \cref{def:direction_quasi_pseudo_normality}
coincide with directional pseudo- and quasi-normality in case of standard constraint mappings from \cite{BenkoCervinkaHoheisel2019}.

We begin by a general result relying on calmness of the constraint function.
Note that we consider the particular situation $\bar y:=0$ for simplicity of
notation. This is not restrictive since $\Phi$ can be shifted appropriately
if $\bar y$ does not vanish to achieve this setting.
Furthermore, we only focus on the concept of directional quasi-normality in
our subsequently stated analysis. Analogous results can be obtained for
directional pseudo-normality.

\begin{proposition}\label{pro:quasi_normality_for_constraint_maps_calm_case}
  A constraint mapping  $\Phi\colon\mathbb X\tto\mathbb Y$ given by $\Phi(x) := g(x) - D$ 
  for all $x\in\mathbb X$,
  where $g\colon\mathbb X \to \mathbb Y$ is calm in direction $u \in \mathbb S_{\mathbb X} $
  at $\bar x\in\mathbb X$ such that $(\bar x, 0) \in \gph \Phi$ 
  and $D \subset \mathbb Y$ is closed,
  is quasi-normal in direction $u$ at $(\bar x, 0)$
  w.r.t.\ some orthonormal basis $\mathcal E:=\{e_1,\ldots,e_m\}\subset\mathbb Y$ of $\mathbb Y$
  provided there do not exist a direction $v \in \mathbb Y$ and 
  a nonzero vector $\lambda\in \mathcal N_D(g(\bar x);v)$ 
  with $0 \in D^*g(\bar x;(u,v))(\lambda)$ satisfying the following condition: 
  there are sequences $\{x_k\}_{k\in\N}\subset\mathbb X$
  with $x_k\neq\bar x$ for all $k\in\N$, $\{z_k\}_{k\in\N}\subset D$,
  $\{\lambda_k\}_{k\in\N}\subset\mathbb Y$, and $\{\eta_k\}_{k\in\N}\subset\mathbb X$
  satisfying $x_k\to\bar x$, $z_k\to g(\bar x)$, $\lambda_k\to\lambda$, $\eta_k\to 0$,
	\begin{equation}\label{eq:directional_convergences_quasi_normality_2}
		\frac{x_k-\bar x}{\norm{x_k-\bar x}}\to u,
		\qquad
		\frac{z_k-g(\bar x)}{\norm{x_k-\bar x}}\to v,
                \qquad
                \frac{g(x_k)-g(\bar x)}{\norm{x_k-\bar x}}\to v,
	\end{equation}
  and, for all $k\in\N$ and $i\in\{1,\ldots,m\}$, $\eta_k \in \widehat{D}^*g(x_k)(\lambda_k)$,
  $\lambda_k\in\widehat{\mathcal N}_D(z_k)$, 
  as well as $\dual{\lambda}{e_i}\dual{g(x_k)-z_k}{e_i}>0$
  if $\dual{\lambda}{e_i}\neq 0$.
  
  Moreover, if $g$ is even calm (particularly Lipschitz continuous) near $\bar x$,
  the two conditions are equivalent.
\end{proposition}
\begin{proof}
    $[\Longleftarrow]$ Choose $\lambda\in\ker D^*\Phi((\bar x,0);(u,0))$ 
    and sequences $\{(x_k,y_k)\}_{k\in\N}\subset\gph\Phi$ with $x_k\neq\bar x$ 
    for all $k\in\N$ and
	$\{\lambda_k\}_{k\in\N}\subset\mathbb Y$, $\{\eta_k\}_{k\in\N}\subset\mathbb X$ satisfying
	\eqref{eq:convergences_definition_quasi_normality} with $\bar y:=0$ and, 
	for all $k\in\N$ and $i\in\{1,\ldots,m\}$,
	$\eta_k\in \widehat{D}^*\Phi(x_k,y_k)(\lambda_k)$ 
	as well as $\dual{\lambda}{e_i}\dual{y_{k}}{e_i}>0$ if $\dual{\lambda}{e_i}\neq 0$.
	Applying 
	\cref{lem:coderivatives_constraint_maps}\,\ref{item:constraint_maps_regular_coderivative}
	yields $\eta_k \in \widehat{D}^*g(x_k)(\lambda_k)$ and
	$\lambda_k\in\widehat{\mathcal N}_D(g(x_k) - y_k)$ for each $k\in\N$.
	The assumed calmness of $g$ at $\bar x$ in direction $u$ yields boundedness
	of the sequence $\{(g(x_k) - g(\bar x))/\norm{x_k - \bar x}\}_{k\in\N}$,
	i.e., along a subsequence (without relabeling) it converges to some $v\in\mathbb Y$.
	Note also that
	$(u,v)\in\mathcal T_{\gph g}(\bar x,g(\bar x))$, i.e., $v\in Dg(\bar x)(u)$,
	and that $\{(x_k,g(x_k))\}_{k\in\N}$ converges to $(\bar x,g(\bar x))$ 
	from direction $(u,v)$.
	Setting $z_k := g(x_k) - y_k$ for each $k\in\N$, 
	we get $z_k \to g(\bar x)$ by continuity of $g$
	as well as $\lambda_k\in\widehat{\mathcal N}_D(z_k)$ 
	and $\dual{\lambda}{e_i}\dual{g(x_k)-z_k}{e_i}>0$ if $\dual{\lambda}{e_i}\neq 0$
	for each $k\in\N$ and $i\in\{1,\ldots,m\}$.
	Moreover, we have
	\begin{align*}
  		\frac{z_k - g(\bar x)}{\norm{x_k-\bar x}}
  		=
  		\frac{g(x_k)- g(\bar x)}{\norm{x_k-\bar x}}
  		-
  		\frac{y_k}{\norm{x_k-\bar x}}
  		\to 
  		v - 0 
  		= 
  		v
  	\end{align*}
  	and $v\in \mathcal T_D(g(\bar x))$ follows as well.
    Finally, taking the limit yields
    $\lambda\in\mathcal N_D(g(\bar x);v)$ and $0 \in D^*g(\bar x;(u,v))(\lambda)$,
    so that the assumptions of the proposition imply $\lambda=0$.
    Consequently, $\Phi$ is quasi-normal in direction $u$ at $(\bar x,0)$ w.r.t.\ $\mathcal E$.
    \\
	$[\Longrightarrow]$ Assume that quasi-normality in direction $u$ holds at 
	$(\bar x,0)$ w.r.t.\ $\mathcal E$ and that $g$ is calm around $\bar x$.
	Suppose that there are some $v \in \mathbb Y$,
	$\lambda\in\mathcal N_D(g(\bar x);v)$ with $0 \in D^*g(\bar x;(u,v))(\lambda)$,
	and sequences $\{x_k\}_{k\in\N}\subset\mathbb X$ with $x_k\neq\bar x$ 
	for all $k\in\N$ and $\{z_k\}_{k\in\N}\subset D$, 
	$\{\lambda_k\}_{k\in\N}\subset\mathbb Y$, $\{\eta_k\}_{k\in\N}\subset\mathbb X$ with
	$x_k\to\bar x$, $z_k\to g(\bar x)$, $\lambda_k\to\lambda$, $\eta_k \to 0$,
	\eqref{eq:directional_convergences_quasi_normality_2}, and, for 
	all $k\in\N$ and $i\in\{1,\ldots,m\}$, $\eta_k \in \widehat{D}^*g(x_k)(\lambda_k)$,
	$\lambda_k\in\widehat{\mathcal N}_D(z_k)$, as well as 
	$\dual{\lambda}{e_i}\dual{g(x_k)-z_k}{e_i}>0$
  	as soon as $\dual{\lambda}{e_i}\neq 0$.
  	Set $y_k:=g(x_k)-z_k$ for each $k\in\N$. Then we have $y_k\to 0$,
  	\begin{align*}
  		\frac{y_k}{\norm{x_k-\bar x}}
  		=
  		\frac{g(x_k)-z_k}{\norm{x_k-\bar x}}
  		=
  		\frac{g(x_k)-g(\bar x)}{\norm{x_k-\bar x}}
  		-\frac{z_k-g(\bar x)}{\norm{x_k-\bar x}}
  		\to v - v = 0,
  	\end{align*}
  	and, for all $k\in\N$ and $i\in\{1,\ldots,m\}$, 
  	$\lambda_k\in\widehat{\mathcal N}_D(g(x_k)-y_k)$ 
  	as well as $\dual{\lambda}{e_i}\dual{y_{k}}{e_i}>0$ if $\dual{\lambda}{e_i}\neq 0$.
  	Since $\eta_k \in \widehat{D}^*g(x_k)(\lambda_k)$,
  	calmness of $g$ at $x_k$ implies $\eta_k\in \widehat{D}^*\Phi(x_k,y_k)(\lambda_k)$
  	due to
  	\cref{lem:coderivatives_constraint_maps}\,\ref{item:constraint_maps_regular_coderivative},
  	and taking the limit yields $\lambda\in\ker D^*\Phi((\bar x,\bar y);(u,0))$.
  	Thus, the assumed quasi-normality of $\Phi$ at $(\bar x,0)$ 
  	in direction $u$ w.r.t.\ $\mathcal E$ yields $\lambda=0$ and the claim follows.
\end{proof}

If $g$ is continuously differentiable, the situation becomes a bit simpler
and we precisely recover the notion of directional quasi-normality for geometric constraint
systems as discussed in \cite[Definition~3.4]{BenkoCervinkaHoheisel2019}.

\begin{corollary}\label{cor:quasi_normality_for_constraint_maps_smooth_case}
  A constraint mapping $\Phi\colon\mathbb X\tto\mathbb Y$ given by $\Phi(x) = g(x) - D$ 
  for all $x\in\mathbb X$, 
  where $g\colon\mathbb X \to \mathbb Y$ is continuously differentiable and 
  $D \subset \mathbb Y$ is closed,
  is quasi-normal in direction $u \in \mathbb S_{\mathbb X}$ at $(\bar x, 0) \in \gph \Phi$ 
  w.r.t.\ some orthonormal basis $\{e_1,\ldots,e_m\}\subset\mathbb Y$ of $\mathbb Y$
  if and only if there does not exist a nonzero vector 
  $\lambda\in\mathcal N_D(g(\bar x);\nabla g(\bar x) u)$ 
  with $ \nabla g(\bar x)^*\lambda=0$
  satisfying the following condition: 
  there are sequences $\{x_k\}_{k\in\N}\subset\mathbb X$ with $x_k\neq\bar x$ 
  for all $k\in\N$, $\{z_k\}_{k\in\N}\subset D$, and $\{\lambda_k\}_{k\in\N}\subset\mathbb Y$
  satisfying $x_k\to\bar x$, $z_k\to g(\bar x)$, $\lambda_k\to\lambda$, 
	\begin{equation}\label{eq:directional_convergences_quasi_normality_3}
		\frac{x_k-\bar x}{\norm{x_k-\bar x}}\to u,
		\qquad
		\frac{z_k-g(\bar x)}{\norm{x_k-\bar x}}\to \nabla g(\bar x)u,
	\end{equation}
  and, for all $k\in\N$ and $i\in\{1,\ldots,m\}$, 
  $\lambda_k\in\widehat{\mathcal N}_D(z_k)$ as well as 
  $\dual{\lambda}{e_i}\dual{g(x_k)-z_k}{e_i}>0$
  if $\dual{\lambda}{e_i}\neq 0$.
\end{corollary}

In \cite[Section~3.3]{BenkoCervinkaHoheisel2019}, it has been reported that under additional
conditions on the set $D$, we can drop the sequences $\{z_k\}_{k\in\N}$ and
$\{\lambda_k\}_{k\in\N}$ from the characterization of directional
quasi-normality in \cref{cor:quasi_normality_for_constraint_maps_smooth_case}.
Particularly, this can be done for so-called \emph{ortho-disjunctive}
programs which cover e.g.\ standard nonlinear, complementarity-, cardinality-, or
switching-constrained optimization problems. 
In this regard, \cref{cor:quasi_normality_for_constraint_maps_smooth_case}
reveals that some results from 
\cite{BertsekasOzdaglar2002,Hestenes1975,KanzowRaharjaSchwartz2021b,KanzowSchwartz2010,LiangYe2021}
are covered by our general concept from \cref{def:direction_quasi_pseudo_normality}.

Let us briefly compare our results with the approach from \cite{BaiYeZhang2019}.

\begin{remark}
	Let us consider the setting discussed in \cref{pro:quasi_normality_for_constraint_maps_calm_case}.
    The directional versions of quasi- and pseudo-normality from \cite{BaiYeZhang2019} operate with all nonzero pairs
    of directions $(u,v)$, rather than just a fixed $u$.
    The advantage is that calmness of $g$ plays no role.
    The reason is, however, that the authors in \cite{BaiYeZhang2019} only derive statements
    regarding metric subregularity, but not metric
    subregularity in some fixed direction.
    Calmness of $g$ is needed precisely for preservation
    of directional information.
    We believe that it is useful to know how to verify if a mapping
    is subregular in a specific direction since only
    some directions play a role in many situations.
    We could drop the calmness assumption from \cref{pro:quasi_normality_for_constraint_maps_calm_case}, but, 
    similarly as in \cite[Theorem~3.1]{BenkoGfrererOutrata2019}, additional directions of the type $(0,v)$ for a nonzero $v$ would appear.
    Clearly, such directions are included among all nonzero pairs $(u,v)$, but the connection to the original direction $u$ 
    would have been lost.
    
    Let us mention that some of the comments from \cite{BaiYeZhang2019}
    about improving \cite[Proposition~2.2]{BenkoGfrererOutrata2019} are not accurate
    since these results are actually not comparable.
    Moreover, e.g. \cite[Corollary~3.1]{BaiYeZhang2019}
    can be easily derived on the basis of \cite[Theorem~3.1]{BenkoGfrererOutrata2019}.
\end{remark}

\subsection{Pseudo- and quasi-normality for problems with equilibrium constraints}

In mathematical optimization, problems with so-called \emph{equilibrium} constraints are used 
to model situations where some variables need to solve a given variational problem.
Exemplary, this covers optimization problems with variational inequality constraints,
see e.g.\ \cite{FacchneiPang2003,LuoPangRalph1996,OutrataKocvaraZowe1998}, or 
bilevel optimization problems, see e.g.\
\cite{Dempe2002,DempeKalashnikovPerezValdesKalashnykova2015}.
In order to model such problems, we need to split the decision space into two parts, i.e.,
we assume that $\mathbb X=\mathbb X^1\times\mathbb X^2$ for Euclidean spaces 
$\mathbb X^1$, $\mathbb X^2$ and exploit $x:=(x^1,x^2)$ for $x\in\mathbb X$, $x^1\in\mathbb X^1$,
and $x^2\in\mathbb X^2$. Furthermore, let $S\colon\mathbb X^1\tto\mathbb X^2$ be the
solution mapping of the underlying variational problem and assume that $\gph S$ is closed. 
For some locally Lipschitz continuous function
$\varphi\colon\mathbb X\to\R$ and some closed set $\Omega\subset\mathbb X^1$, 
the problem of interest is given by
\begin{equation}\label{eq:MPEC}\tag{MPEC}
	\min\{\varphi(x)\,|\,x^1\in\Omega,\,x^2\in S(x^1)\}.
\end{equation}
Introducing $\Phi\colon\mathbb X\tto\mathbb X$ by means of
\begin{equation}\label{eq:Phi_for_MPEC}
	\forall x\in\mathbb X\colon\quad
	\Phi(x):=\bigl(\Omega-x^1,S(x^1)-x^2\bigr),
\end{equation}
we easily see that \eqref{eq:MPEC} is a special instance of \eqref{eq:nonsmooth_problem}
with $\bar y:=0$.

In order to apply our new notions of directional pseudo- and quasi-normality 
from \cref{def:direction_quasi_pseudo_normality} to
\eqref{eq:MPEC}, we need to compute the regular and directional limiting coderivative
of $\Phi$ from \eqref{eq:Phi_for_MPEC}.
For $(x,y)\in\gph\Phi$ and $(\lambda,\mu)\in\mathbb X$, 
\cite[Lemma~5.7(ii)]{BenkoMehlitz2020} yields
\begin{align*}
	\widehat D^*\Phi(x,y)(\lambda,\mu)
	=
	\begin{cases}
		\bigl(\widehat D^*S(x^1,x^2+y^2)(\mu)-\lambda,-\mu\bigr)	&	-\lambda\in\widehat{\mathcal N}_\Omega(x^1+y^1),\\
		\varnothing												& \text{otherwise.}
	\end{cases}
\end{align*}
Furthermore, for $(x,0)\in\gph\Phi$, $(\lambda,\mu)\in\mathbb X$, and some direction $u\in\mathbb X$, 
\cite[Lemma~5.6(ii)]{BenkoMehlitz2020} gives us the upper estimate
\begin{align*}
	D^*\Phi((x,0);(u,0))(\lambda,\mu)
	\subset
	\begin{cases}
		\bigl(D^*S(x;u)(\mu)-\lambda,-\mu\bigr)	&	-\lambda\in \mathcal N_\Omega(x^1;u^1),\\
		\varnothing												& \text{otherwise.}
	\end{cases}
\end{align*}
Particularly, we obtain 
\[
	\ker D^*\Phi((x,0);(u,0))
	\subset
	\bigl(D^*S(x;u)(0)\cap\bigl(-\mathcal N_\Omega(x^1;u^1)\bigr)\bigr)\times\{0\}.
\]
With this in mind, we have the following result.
\begin{lemma}\label{lem:pseudo_normality_for_MPEC}
	Fix $(\bar x,0)\in\gph\Phi$ where $\Phi$ is given in \eqref{eq:Phi_for_MPEC}
	as well as a direction $u\in\mathbb S_{\mathbb X}$.
	Then the following statements hold.
	\begin{enumerate}
		\item\label{item:pseudo_normality_MPEC}
			Suppose that there does not exist
			a nonzero $\lambda\in D^*S(x;u)(0)\cap\bigl(-\mathcal N_\Omega(x^1;u^1)\bigr)$
			satisfying the following condition: 
			there are sequences $\{(x_k,y_k)\}_{k\in\N}\subset\gph\Phi$ with $x_k\neq\bar x$ for all $k\in\N$
			and $\{\lambda_k\}_{k\in\N},\{\eta_k\}_{k\in\N}\subset\mathbb X^1$, $\{\mu_k\}_{k\in\N}\subset\mathbb X^2$,
			such that \eqref{eq:convergences_definition_quasi_normality} and $\mu_k\to 0$
			as well as $\eta_k+\lambda_k\in\widehat D^*S(x_k^1,x_k^2+y_k^2)(\mu_k)$, $-\lambda_k\in\widehat{\mathcal N}_\Omega(x_k^1+y_k^1)$,
			and $\dual{\lambda}{y_k^1}>0$ for all $k\in\N$.
			Then $\Phi$ is pseudo-normal at $(\bar x,0)$ in direction $u$.
		\item Let $\mathcal E^j:=\{e^j_1,\ldots,e^j_{m_j}\}$ be an orthonormal basis of $\mathbb X^j$ for $j=1,2$.
			Suppose that there does not exist 
			a nonzero $\lambda\in D^*S(x;u)(0)\cap\bigl(-\mathcal N_\Omega(x^1;u^1)\bigr)$
			satisfying the following condition: 
			there are sequences $\{(x_k,y_k)\}_{k\in\N}\subset\gph\Phi$ with $x_k\neq\bar x$ for all $k\in\N$
			and $\{\lambda_k\}_{k\in\N},\{\eta_k\}_{k\in\N}\subset\mathbb X^1$, $\{\mu_k\}_{k\in\N}\subset\mathbb X^2$,
			such that \eqref{eq:convergences_definition_quasi_normality} and $\mu_k\to 0$
			as well as, for all $k\in\N$ and $i\in\{1,\ldots,m_1\}$,
			$\eta_k+\lambda_k\in\widehat D^*S(x_k^1,x_k^2+y_k^2)(\mu_k)$, $-\lambda_k\in\widehat{\mathcal N}_\Omega(x_k^1+y_k^1)$,
			and $\dual{\lambda}{e^1_i}\dual{y_k^1}{e^1_i}>0$ if $\dual{\lambda}{e^1_i}\neq 0$.
			Then $\Phi$ is quasi-normal at $(\bar x,0)$ in direction $u$ w.r.t.\ the orthonormal basis $\mathcal E^1\times\mathcal E^2$
			of $\mathbb X$.
	\end{enumerate}
\end{lemma}

We note that, depending on the structure of the underlying variational problem, the appearing coderivatives of $S$ can be
specified or at least estimated from above in terms of initial problem data under mild assumptions, see e.g.\
\cite{DontchevRockafellar2014,GfrererOutrata2016,GfrererOutrata2016c,Mordukhovich2006} and the references therein.
Observe that working with upper estimates of these derivatives in \cref{lem:pseudo_normality_for_MPEC} still
yields sufficient conditions for directional pseudo- and quasi-normality of $\Phi$.
Clearly, when applied to bilevel optimization problems, \cref{lem:pseudo_normality_for_MPEC} provides suitable
notions of directional pseudo- and quasi-normality for this problem class in hierarchical form.
Combined with \cref{lem:directional_M_stationarity_via_metric_subregularity} and
\cref{thm:quasi_normality_yields_asymptotic_regularity}, validity of the conditions from 
\cref{lem:pseudo_normality_for_MPEC} guarantees directional M-stationarity of local minimizers
associated with \eqref{eq:MPEC}.
Let us briefly mention \cite{BaiYe2021} where another notion of directional quasi-normality has been 
introduced for bilevel optimization problems which is based on the so-called
value function reformulation. The latter is a single-level optimization problem with highly irregular
nonsmooth inequality constraints and essentially different from \eqref{eq:MPEC}.
However, a common feature of both approaches is that the considered problems are only implicitly
given. While in \eqref{eq:MPEC}, $S$ is an implicit object, the same holds true for the value
function. Using tools of generalized differentiation allows us to characterize the derivatives
of these objects and, thus, end up with explicit conditions in terms of problem data.

Let us investigate a small example to illustrate the conditions from \cref{lem:pseudo_normality_for_MPEC}.
\begin{example}\label{ex:pseudo_normality_for_MPEC}
	Let $\mathbb X^1:=\mathbb X^2:=\R$ and $S\colon\R\tto\R$ be given by
	\[
		\forall t\in\R\colon\quad
		S(t)
		:=
		\begin{cases}
			\{-t^2\}	&	t\leq 0,\\
			\{\sqrt t\}	&	t>0.
		\end{cases}
	\]
	We consider the feasible region
	\[
		\{x\in\R^2\,|\,x^1\in\R_+,\,x^2\in S(x^1)\}
	\]
	at $\bar x:=(0,0)$.
	
	One obtains $\mathcal T_{\gph S}(\bar x)=(\R_-\times\{0\})\cup(\{0\}\times\R_+)$ and $\mathcal T_{\R_+}(\bar x^1)=\R_+$,
	so that the only interesting direction from $\mathbb S_{\R^2}$ is $u:=(0,1)$.
	Note that we have $D^*S(\bar x;u)(0)\cap(-\mathcal N_{\R_+}(\bar x^1,u^1))=\R_+$.
	Suppose that there are $\lambda>0$
	and sequences $\{(x_k,y_k)\}_{k\in\N}\subset\gph\Phi$, $\{\lambda_k\}_{k\in\N},\{\eta_k\}_{k\in\N},\{\mu_k\}_{k\in\N}\subset\R$
	as in statement~\ref{item:pseudo_normality_MPEC} of \cref{lem:pseudo_normality_for_MPEC}.
	Due to $\lambda_k\to\lambda$, $-\lambda_k\in\widehat{\mathcal N}_{\R_+}(x_k^1+y_k^1)$, and $\lambda y_k^1>0$, 
	we find $x_k^1=-y_k^1<0$ for each $k\in\N$.
	Thus, we find $S(x_k^1)=\{-(x_k^1)^2\}$ for each $k\in\N$, and an evaluation of the coderivative condition regarding
	$S$ yields $\eta_k+\lambda_k=-2x_k^1\mu_k$ for each $k\in\N$. Taking the limit $k\to\infty$ yields $\lambda=0$ which
	is a contradiction.
	Thus, due to \cref{lem:pseudo_normality_for_MPEC}\,\ref{item:pseudo_normality_MPEC}, the associated mapping $\Phi$
	from \eqref{eq:Phi_for_MPEC} is pseudo-normal w.r.t.\ all directions from $\mathbb S_{\R^2}$.
\end{example}

\section{Sufficient conditions for asymptotic regularity via pseudo-coderivatives}\label{sec:asymptotic_regularity_via_super_coderivative}

\subsection{On the super-coderivative}\label{sec:super_coderivative}

As we will observe in this section, (strong) directional asymptotic regularity 
can be suitably investigated via
the following novel concept of generalized differentiation, which generalizes the notion of pseudo-coderivatives.
\begin{definition}\label{def:super_coderivative}
	Let $\Phi\colon\mathbb X\tto\mathbb Y$ be a set-valued mapping with a closed graph
	and fix $(\bar x,\bar y)\in\gph\Phi$ and 
	$(u,v)\in\mathbb S_{\mathbb X}\times\mathbb S_{\mathbb Y}$.
	The \emph{super-coderivative} of $\Phi$ at $(\bar x,\bar y)$ in direction $(u,v)$ is
	the mapping $D^*_\textup{sup}\Phi((\xb,\yb); (u,v))\colon\mathbb Y \tto \mathbb X$,
	which assigns to every $y^* \in\mathbb Y$ the set of all $x^\ast \in \mathbb X$ 
	for which there exist sequences
	$\{u_k\}_{k\in\N},\{x_k^*\}_{k\in\N}\subset\mathbb X$, 
	$\{v_k\}_{k\in\N},\{y_k^*\}_{k\in\N}\subset\mathbb Y$,
	and $\{t_k\}_{k\in\N},\{\tau_k\}_{k\in\N}\subset\R_+$ which satisfy
	$u_k\to u$, $v_k\to v$, $x_k^*\to x^*$, $y_k^*\to y^*$, $t_k\searrow 0$, 
	$\tau_k\searrow 0$, and $\tau_k/t_k\to 0$
	such that
	\begin{equation*}
		(\tau_k \norm{v_k}) x_k^* 
		\in 
		\widehat{D}^*\Phi(\xb + t_k u_k,\yb + \tau_k v_k)((t_k \norm{u_k}) y_k^*)
	\end{equation*}
	or, equivalently,
	\begin{equation}\label{eq:characterization_super_coderivative}
		x_k^* 
		\in 
		\widehat{D}^*\Phi(\xb + t_k u_k,\yb + \tau_k v_k)
			(((t_k \norm{u_k})/(\tau_k \norm{v_k})) y_k^*)
	\end{equation}
	holds for all $k\in\N$.
\end{definition}

We start with some remarks regarding \cref{def:super_coderivative}.
First, observe that we only exploit the super-coderivative w.r.t.\
unit directions $(u,v)\in\mathbb S_{\mathbb X}\times\mathbb S_{\mathbb Y}$
which also means that $\{u_k\}_{k\in\N}\subset\mathbb X$ and
$\{v_k\}_{k\in\N}\subset\mathbb Y$ can be chosen such that $u_k\neq 0$
and $v_k\neq 0$ hold for all $k\in\N$. 
Particularly, condition \eqref{eq:characterization_super_coderivative} is reasonable.

Second, we would like to note that 
$x^*\in D^*_\textup{sup}\Phi((\bar x,\bar y);(u,v))(y^*)$ implies the existence of
sequences $\{u_k\}_{k\in\N}\subset\mathbb X$, 
$\{v_k\}_{k\in\N}\subset\mathbb Y$,
and $\{t_k\}_{k\in\N},\{\tau_k\}_{k\in\N}\subset\R_+$ which satisfy
$u_k\to u$, $v_k\to v$, $t_k\searrow 0$, $\tau_k\searrow 0$, and $\tau_k/t_k\to 0$
as well as $(\bar x+t_ku_k,\bar y+\tau_kv_k)\in\gph\Phi$ for all $k\in\N$.
Thus, in the light of the definition of the graphical subderivative, one might be tempted
to say that the pair $(u,v)$ belongs to the graph of the graphical
\emph{super-derivative} of $\Phi$ at $(\bar x,\bar y)$.
This justifies the terminology in \cref{def:super_coderivative}.

Let us briefly discuss the relation between pseudo-coderivatives and the
novel super-coderivative from \cref{def:super_coderivative}.
Consider $\gamma > 1$ and $x^* \in D^*_\gamma\Phi((\xb,\yb); (u,v))(y^*)$ 
for $(u,v)\in\mathbb S_{\mathbb X}\times\mathbb S_{\mathbb Y}$, and $y^*\in\mathbb Y^*$.
Setting $\tau_k := (t_k \norm{u_k})^{\gamma}$ for each $k\in\N$, 
where $\{t_k\}_{k\in\N}\subset\R_+$ 
and $\{u_k\}_{k\in\N}\subset\mathbb X$ are the sequences from the definition of the 
pseudo-coderivative,
we get $x^* \in D^*_\textup{sup}\Phi((\xb,\yb); (u,v))(y^*)$ 
since $t_k^{\gamma-1}\norm{u_k}^\gamma\to 0$.

In the subsequent lemma, we comment on the converse inclusion which, to some extent, 
holds in the presence of a qualification condition.
\begin{lemma}\label{lem:super_coderivative_vs_pseudo_coderivative}
	Let $(\bar x,\bar y)\in\gph\Phi$, $(u,v)\in\mathbb S_{\mathbb X}\times\mathbb S_{\mathbb Y}$,
	$y^* \in \mathbb Y$, and $\gamma > 1$ be fixed.
	Furthermore, assume that $\ker D^*_\gamma\Phi((\bar x,\bar y);(u,0))=\{0\}$ holds.
	Then there exists $\alpha > 0$ such that
	\begin{align*}
		D^*_\textup{sup}\Phi((\bar x,\bar y);(u,v))(y^*)
		&\subset
		\widetilde{D}^*_\gamma\Phi((\bar x,\bar y);(u,0))(0)
		\cup
		D^*_\gamma\Phi((\bar x,\bar y);(u,\alpha v))(y^*/\alpha)
		\\
		&\qquad
		\cup
		\Im D^*_\gamma\Phi((\bar x,\bar y);(u,0))
		\\
		&\subset
		\Im \widetilde D^*_\gamma\Phi((\bar x,\bar y);(u,0)).
	\end{align*}
\end{lemma}
\begin{proof}
	Let $x^*\in D^*_\textup{sup}\Phi((\bar x,\bar y);(u,v))(y^*)$ be arbitrarily chosen.
	Then we find sequences $\{u_k\}_{k\in\N},\{x_k^*\}_{k\in\N}\subset\mathbb X$, 
	$\{v_k\}_{k\in\N},\{y_k^*\}_{k\in\N}\subset\mathbb Y$,
	and $\{t_k\}_{k\in\N},\{\tau_k\}_{k\in\N}\subset\R_+$ which satisfy
	$u_k\to u$, $v_k\to v$, $x_k^*\to x^*$, $y_k^*\to y^*$, $t_k\searrow 0$, 
	$\tau_k\searrow 0$, and $\tau_k/t_k\to 0$
	as well as \eqref{eq:characterization_super_coderivative} for all $k\in\N$.
	This also gives us
	\begin{equation}\label{eq:insert_order_gamma_into_super_coderivative}
		x_k^*
		\in
		\widehat D^*\Phi
			\left(
				\bar x+t_ku_k,
				\bar y+(t_k\norm{u_k})^\gamma\frac{\tau_kv_k}{(t_k\norm{u_k})^\gamma}
			\right)
			\left(
				(t_k\norm{u_k})^{1-\gamma}
				\frac{(t_k\norm{u_k})^\gamma}{\tau_k\norm{v_k}}y_k^*
			\right)
	\end{equation}
	for all $k\in\N$.
	Set $\tilde y_k^*:=(t_k\norm{u_k})^\gamma/(\tau_k\norm{v_k})y_k^*$ for each $k\in\N$. 
	In case where $\{\tilde y_k^*\}_{k\in\N}$ is not bounded, 
	we have $(\tau_k\norm{v_k})/(t_k\norm{u_k})^\gamma\to 0$
	along a subsequence (without relabeling), and taking the limit in
	\[
		x_k^*/\nnorm{\tilde y_k^*}
		\in
		\widehat D^*\Phi
		\left(
			\bar x+t_ku_k,
			\bar y+(t_k\norm{u_k})^\gamma\frac{\tau_kv_k}{(t_k\norm{u_k})^\gamma}
		\right)
		\left(
			(t_k\norm{u_k})^{1-\gamma} 
			\tilde y_k^*/\nnorm{\tilde y_k^*}
		\right)
	\]
	yields $\ker D^*_\gamma\Phi((\bar x,\bar y);(u,0))\neq\{0\}$ which is a contradiction.
	Hence, $\{\tilde y_k^*\}_{k\in\N}$ is bounded.
	
	For each $k\in\N$, we set $\alpha_k:=\tau_k\norm{v_k}/(t_k\norm{u_k})^\gamma$.
	First, suppose that $\{\alpha_k\}_{k\in\N}$ is not bounded.
	Then, along a subsequence (without relabeling), we may assume $\alpha_k\to\infty$.
	By boundedness of $\{y_k^*\}_{k\in\N}$, $\tilde y_k^*\to 0$ follows.
	Rewriting \eqref{eq:insert_order_gamma_into_super_coderivative} yields
	\[
		x_k^*\in\widehat D^*\Phi\left(\bar x+t_ku_k,\bar y+t_k\frac{\tau_kv_k}{t_k}\right)
		\left((t_k\norm{u_k})^{1-\gamma}\tilde y_k^*\right)
	\]
	for each $k\in\N$, and taking the limit $k\to\infty$ 
	while respecting $\tau_k/t_k\to 0$, thus, gives 
	$x^*\in \widetilde D^*_\gamma\Phi((\bar x,\bar y);(u,0))(0)$.
	In case where $\{\alpha_k\}_{k\in\N}$ converges to some $\alpha>0$ (along a subsequence
	without relabeling), we can simply take the limit $k\to\infty$ in
	\eqref{eq:insert_order_gamma_into_super_coderivative} in order to find
	$x^*\in D^*_\gamma\Phi((\bar x,\bar y);(u,\alpha v))(y^*/\alpha)$.
	Finally, let us consider the case $\alpha_k\to 0$ (along a subsequence without
	relabeling). Then, by boundedness of $\{\tilde y_k^*\}_{k\in\N}$,
	taking the limit $k\to\infty$ in \eqref{eq:insert_order_gamma_into_super_coderivative}
	gives $x^*\in \Im D^*_\gamma\Phi((\bar x,\bar y);(u,0))$.
	Thus, we have shown the first inclusion.
	
	The second inclusion follows by the trivial upper estimate for the pseudo-coderivative.
\end{proof}

Let us now interrelate the concepts of super-coderivatives and asymptotic regularity.
Choose sequences $\{(x_k,y_k)\}_{k\in\N}\subset\gph\Phi$, 
$\{x_k^*\}_{k\in\N}\subset\mathbb X$, 	and $\{\lambda_k\}_{k\in\N}\subset\mathbb Y$ 
as well as $x^*\in\mathbb X$ and $y^*\in\mathbb Y$ satisfying
$x_k\notin\Phi^{-1}(\bar y)$, $y_k\neq\bar y$, 
and $x_k^*\in\widehat{D}^*\Phi(x_k,y_k)(\lambda_k)$
for all $k\in\N$ as well as the convergences
\eqref{eq:convergences_directional_asymptotic_regularity}.
For each $k\in\N$, we set $t_k:=\norm{x_k-\bar x}$, $\tau_k:=\norm{y_k-\bar y}$,
\[
	u_k:=\frac{x_k-\bar x}{\norm{x_k-\bar x}},\qquad 
	v_k:=\frac{y_k-\bar y}{\norm{y_k-\bar y}},\qquad
	y_k^*:=\frac{\norm{y_k-\bar y}}{\norm{x_k-\bar x}}\lambda_k
\]
and find
\[
	\forall k\in\N\colon\quad
		x_k^*\in \widehat{D}^*\Phi(\bar x+t_ku_k,\bar y+\tau_kv_k)((t_k/\tau_k)y_k^*).
\]
Along a subsequence (without relabeling), $v_k\to v$ holds for some 
$v\in\mathbb S_{\mathbb Y}$.
Thus, taking the limit $k\to\infty$, 
we have $x^*\in D^*_\textup{sup}\Phi((\bar x,\bar y);(u,v))(y^*)$
by definition of the super-coderivative.
Moreover, from \eqref{eq:convergences_directional_asymptotic_regularity} 
we also know that $y^* = \norm{y^*} v$.
Consequently, we come up with the following lemma.

\begin{lemma}\label{lem:super_coderivative_vs_asymptotic_regularity}
	Let $(\bar x,\bar y)\in\gph\Phi$ and $u\in\mathbb S_{\mathbb X}$ be fixed.
	If 
	\[
		\bigcup_{v\in\mathbb S_{\mathbb Y}}
		D^*_\textup{sup}\Phi((\bar x,\bar y);(u,v))(\beta v)
		\subset
		\Im D^*\Phi(\bar x,\bar y) 
	\]
	holds for all $\beta \geq 0$, then $\Phi$ is asymptotically regular at $(\bar x,\bar y)$ in direction $u$.
	If the above estimate holds for all $\beta\geq 0$
	with $\Im D^*\Phi(\bar x,\bar y)$ replaced by $\Im D^*\Phi((\bar x,\bar y);(u,0))$,
	then $\Phi$ is strongly asymptotically regular at $(\bar x,\bar y)$ in direction $u$.
\end{lemma}

The next result follows as a corollary of 
\cref{lem:super_coderivative_vs_asymptotic_regularity,lem:super_coderivative_vs_pseudo_coderivative},
and gives new sufficient conditions for directional asymptotic regularity.
Note that strong directional asymptotic regularity can be handled analogously.
\begin{theorem}\label{thm:asymptotic_regularity_via_pseudo_coderivatives}
	Let $(\bar x,\bar y)\in\gph\Phi$, $u\in\mathbb S_{\mathbb X}$, and $\gamma>1$ be fixed.
	Furthermore, assume that $\ker D^*_\gamma\Phi((\bar x,\bar y);(u,0))=\{0\}$ holds.
	If
	\begin{equation}\label{eq:sufficient_condition_asym_reg_pseudo}
		\widetilde{D}^*_\gamma\Phi((\bar x,\bar y);(u,0))(0)
		\cup
		\bigcup_{v\in\mathbb S_{\mathbb Y}}
		D^*_\gamma\Phi((\bar x,\bar y);(u,\alpha v))(\beta v)
		\subset
		\Im D^*\Phi(\bar x,\bar y)
	\end{equation}
	holds for all $\alpha,\beta \geq 0$, particularly, if
	\begin{equation}\label{eq:sufficient_condition_asym_reg_pseudo_rough}
		\Im\widetilde D^*_\gamma\Phi((\bar x,\bar y);(u,0))
		\subset
		\Im D^*\Phi(\bar x,\bar y)
	\end{equation}
	holds, then $\Phi$ is asymptotically regular at $(\bar x,\bar y)$ in direction $u$.
\end{theorem}

In case where the pseudo-coderivatives involved in the statement of
\cref{thm:asymptotic_regularity_via_pseudo_coderivatives} can be computed
or estimated from above, new applicable sufficient conditions for (strong)
directional asymptotic regularity are at hand.
Particularly, in situations where $\Phi$ is given in form of a constraint
mapping and $\gamma:=2$ is fixed, we can rely on the results obtained
in \cite[Section~3]{BenkoMehlitz2022} in order to make the findings
of \cref{thm:asymptotic_regularity_via_pseudo_coderivatives} more specific.
This will be done in the next subsection.
	
\subsection{Constraint mappings}\label{sec:super_coderivative_condition_constraint_maps}

Throughout the section, we assume that $\Phi\colon\mathbb X\tto\mathbb Y$
is given by $\Phi(x):=g(x)-D$, $x\in\mathbb X$, where $g\colon\mathbb X\to\mathbb Y$ is a
single-valued, twice continuously differentiable function and $D\subset\mathbb Y$
is a closed set. Furthermore, for simplicity of notation, we fix $\bar y:=0$ which
is not restrictive as already mentioned earlier.

We start with a general result which does not rely on any additional structure
of the set $D$.
\begin{theorem}\label{thm:asymptotic_regularity_pseudo_coderivative_non_polyhedral}
	Let $(\bar x,0)\in\gph\Phi$ as well as $u\in\mathbb S_{\mathbb X}$
	be fixed. Assume that the condition
	\begin{equation}\label{eq:CQ_pseudo_subregularity_nonpolyhedral_II}
		\left.
		\begin{aligned}
			&\nabla g(\bar x)^*y^*=0,\,
			\nabla^2\dual{y^*}{g}(\bar x)(u)+\nabla g(\bar x)^* z^*=0,\\
			& y^*\in\mathcal N_D(g(\bar x);\nabla g(\bar x)u),\,
			z^*\in D\mathcal N_D(g(\bar x),y^*)(\nabla g(\bar x)u)
		\end{aligned}
		\right\}
		\quad\Longrightarrow\quad
		y^*=0
	\end{equation}
	is valid. Furthermore, let 
	\begin{equation}\label{eq:CQ_pseudo_subregularity_nonpolyhedral_Ia}
		\left.
		\begin{aligned}
			&\nabla g(\bar x)^*y^*=0,\,
			\nabla g(\bar x)^*\hat z^*=0,\\
			& y^*\in\mathcal N_D(g(\bar x);\nabla g(\bar x)u),\,
			\hat z^*\in D\mathcal N_D(g(\bar x),y^*)(0)
		\end{aligned}
		\right\}
		\quad\Longrightarrow\quad
		\hat z^*=0
	\end{equation}
	or, in case $\nabla g(\bar x)u\neq 0$,
	\begin{equation}\label{eq:CQ_pseudo_subregularity_nonpolyhedral_Ib}
		\left.
		\begin{aligned}
			&\nabla g(\bar x)^*y^*=0,\,
			\nabla g(\bar x)^*\hat z^*=0,\\
			& y^*\in\mathcal N_D(g(\bar x);\nabla g(\bar x)u)
		\end{aligned}
		\right\}
		\quad\Longrightarrow\quad
		\hat z^*\notin D_\textup{sub}\mathcal N_D(g(\bar x),y^*)
			\left(\frac{\nabla g(\bar x)u}{\norm{\nabla g(\bar x)u}}\right)
	\end{equation}
	hold.
	\begin{enumerate}
		\item If, for each $x^*\in\mathbb X$ and $y^*,z^*\in\mathbb Y$ satisfying
			\begin{subequations}\label{eq:system_pseudo_coderivative_nonpolyhedral}
				\begin{align}
					\label{eq:pseudo_coderivative_nonpolyhedral_x}
						x^*
						&=
						\nabla^2\dual{y^*}{g}(\bar x)(u)
						+
						\nabla g(\bar x)^*z^*,
						\\
					\label{eq:pseudo_coderivative_nonpolyhedral_y}
						y^*
						&\in 
						\mathcal N_D(g(\bar x);\nabla g(\bar x)u)
						\cap
						\ker\nabla g(\bar x)^*,
						\\
					\label{eq:pseudo_coderivative_nonpolyhedral_z}
						z^*
						&\in 
						D\mathcal N_D(g(\bar x),y^*)(\nabla g(\bar x)u),
				\end{align}
			\end{subequations}
			there is some $\lambda\in\mathcal N_D(g(\bar x))$
			such that $x^*=\nabla g(\bar x)^*\lambda$,
			then $\Phi$ is asymptotically regular at $(\bar x,0)$ in direction $u$.
		\item If, for each $x^*\in\mathbb X$ and $y^*,z^*\in\mathbb Y$ satisfying
			\eqref{eq:system_pseudo_coderivative_nonpolyhedral}, there is some
			$\lambda\in\mathcal N_D(g(\bar x);\nabla g(\bar x)u)$
			such that $x^*=\nabla g(\bar x)^*\lambda$, then $\Phi$ is
			strongly asymptotically regular at $(\bar x,0)$ in direction $u$.			
	\end{enumerate}
\end{theorem}
\begin{proof}
	Let us start with the proof of the first statement.
	From \cite[Theorem~3.2\,(b)]{BenkoMehlitz2022},
	\eqref{eq:CQ_pseudo_subregularity_nonpolyhedral_II} 
	together with
	\eqref{eq:CQ_pseudo_subregularity_nonpolyhedral_Ia} or, in case $\nabla g(\bar x)u\neq 0$,
	\eqref{eq:CQ_pseudo_subregularity_nonpolyhedral_Ib}
	imply
	$\ker\widetilde D^*_2\Phi((\bar x,0);(u,0))=\{0\}$ which also
	implies $\ker D^*_2\Phi((\bar x,0);(u,0))=\{0\}$.
	Now, pick $x^*\in\Im \widetilde D^*_2\Phi((\bar x,0);(u,0))$.
	Then, due to \cite[Theorem~3.2\,(b)]{BenkoMehlitz2022},
	we find $y^*,z^*\in\mathbb Y$ satisfying 
	\eqref{eq:system_pseudo_coderivative_nonpolyhedral}.
	The assumptions guarantee that we can find $\lambda\in\mathcal N_D(g(\bar x))$
	such that $x^*=\nabla g(\bar x)^*\lambda\in \Im D^*\Phi(\bar x,0)$
	where we used 
	\cref{lem:coderivatives_constraint_maps}\,\ref{item:constraint_maps_limiting_coderivative}.
	It follows $\Im \widetilde D^*_2\Phi((\bar x,0);(u,0))\subset\Im D^*\Phi(\bar x,0)$.
	Thus, \cref{thm:asymptotic_regularity_via_pseudo_coderivatives}
	shows that $\Phi$ is asymptotically regular at $(\bar x,0)$ in direction $u$.
	The second statement follows in analogous way while
	respecting 
	\cref{lem:coderivatives_constraint_maps}\,\ref{item:constraint_maps_directional_limiting_coderivative}.
\end{proof}

We note that \eqref{eq:CQ_pseudo_subregularity_nonpolyhedral_Ia} is stronger than \eqref{eq:CQ_pseudo_subregularity_nonpolyhedral_Ib}
when $\nabla g(\bar x)u\neq 0$ holds, see \cite[formula (2.2)]{BenkoMehlitz2022}.
Naturally, this means that it is sufficient to check \eqref{eq:CQ_pseudo_subregularity_nonpolyhedral_Ia}
regardless whether $\nabla g(\bar x)u$ vanishes or not. In case $\nabla g(\bar x)u\neq 0$, however,
it is already sufficient to check the milder condition \eqref{eq:CQ_pseudo_subregularity_nonpolyhedral_Ib}.
This will be important later on, see \cref{Pro:Milder_than_FOSCMS} and \cref{rem:forgetting_directional_information_in_CQ} below.

The subsequently stated results address the particular case where $D$ is
a polyhedral set, i.e., it is the union of finitely many convex polyhedral sets.
Similarly, $D$ is referred to as locally polyhedral around $y \in D$ whenever 
$D\cap V$ is polyhedral for some neighborhood $V$ of $y$.

\begin{theorem}\label{thm:asymptotic_regularity_pseudo_coderivative_polyhedral_I}
	Let $(\bar x,0)\in\gph\Phi$ as well as $u\in\mathbb S_{\mathbb X}$
	be fixed. Let $\mathbb Y:=\R^m$ and let $D$ be polyhedral
	locally around $g(\bar x)$.
	Assume that the following condition holds:
	\begin{equation}\label{eq:CQ_pseudo_subregularity_polyhedral_I}
		\left.
		\begin{aligned}
			&\nabla g(\bar x)^*y^*=0,\,
			\nabla^2\dual{y^*}{g}(\bar x)(u)+\nabla g(\bar x)^* z^*=0,\\
			& y^*\in\mathcal N_{\mathcal T_D(g(\bar x))}(\nabla g(\bar x)u),\,
			z^*\in \mathcal T_{\mathcal N_{\mathcal T_D(g(\bar x))}(\nabla g(\bar x)u)}(y^*)
		\end{aligned}
		\right\}
		\quad\Longrightarrow\quad
		y^*=0.
	\end{equation}
	\begin{enumerate}
		\item If, for each $x^*\in\mathbb X$ and $y^*,z^*\in\R^m$ satisfying
			\eqref{eq:pseudo_coderivative_nonpolyhedral_x} and
			\begin{equation}\label{eq:system_pseudo_coderivative_polyhedral_I}
				\begin{aligned}
						y^*
						&\in 
						\mathcal N_{\mathcal T_D(g(\bar x))}(\nabla g(\bar x)u)
						\cap
						\ker\nabla g(\bar x)^*,
						\\
						z^*
						&\in 
						\mathcal T_{
							\mathcal N_{\mathcal T_D(g(\bar x))}(\nabla g(\bar x)u)
							}
							(y^*),
				\end{aligned}
			\end{equation}
			there is some $\lambda\in\mathcal N_D(g(\bar x))$
			such that $x^*=\nabla g(\bar x)^*\lambda$,
			then $\Phi$ is asymptotically regular at $(\bar x,0)$ in direction $u$.
		\item If, for each $x^*\in\mathbb X$ and $y^*,z^*\in\R^m$ satisfying
			\eqref{eq:pseudo_coderivative_nonpolyhedral_x} and
			\eqref{eq:system_pseudo_coderivative_polyhedral_I}, there is some
			$\lambda\in\mathcal N_{\mathcal T_D(g(\bar x))}(\nabla g(\bar x)u)$
			such that $x^*=\nabla g(\bar x)^*\lambda$, then $\Phi$ is
			strongly asymptotically regular at $(\bar x,0)$ in direction $u$.			
	\end{enumerate}
\end{theorem}
\begin{proof}
	The proof is analogous to the one of
	\cref{thm:asymptotic_regularity_pseudo_coderivative_non_polyhedral}
	and exploits \cite[Theorem~3.2\,(c)]{BenkoMehlitz2022}.
\end{proof}

In \cref{thm:asymptotic_regularity_pseudo_coderivative_non_polyhedral,thm:asymptotic_regularity_pseudo_coderivative_polyhedral_I},
we relied on the more restrictive assumption \eqref{eq:sufficient_condition_asym_reg_pseudo_rough} from
\cref{thm:asymptotic_regularity_via_pseudo_coderivatives}.
In the general case, we were not able to utilize the more refined condition \eqref{eq:sufficient_condition_asym_reg_pseudo},
but in the polyhedral case, we obtain the following improved result.
We would like to point out that, based on 
\cref{thm:asymptotic_regularity_via_pseudo_coderivatives},
one can state even finer but more technical sufficient conditions
for directional asymptotic regularity in the polyhedral case.

\begin{theorem}\label{thm:asymptotic_regularity_pseudo_coderivative_polyhedral_II}
	Let $(\bar x,0)\in\gph\Phi$ as well as $u\in\mathbb S_{\mathbb X}$
	be fixed. 
	Let $\mathbb Y:=\R^m$ and let $D$ be polyhedral
	locally around $g(\bar x)$.
	Furthermore, we set
	$\mathbf T(u):=\mathcal T_{\mathcal T_D(g(\bar x))}(\nabla g(\bar x)u)$
	and, for arbitrary $s\in\mathbb X$ and $v\in\R^m$,
	$w_s(u,v):=\nabla g(\bar x)s+1/2\nabla^2g(\bar x)[u,u]-v$.
	Assume that the following condition holds for each $s\in\mathbb X$:
	\begin{equation}\label{eq:CQ_pseudo_subregularity_polyhedral_II}
		\left.
		\begin{aligned}
			&\nabla g(\bar x)^*y^*=0,\,
			\nabla^2\dual{y^*}{g}(\bar x)(u)+\nabla g(\bar x)^* z^*=0,\\
			& y^*\in\mathcal N_{\mathbf T(u)}(w_s(u,0)),\,
			z^*\in \mathcal T_{\mathcal N_{\mathbf T(u)}(w_s(u,0))}(y^*)
		\end{aligned}
		\right\}
		\quad\Longrightarrow\quad
		y^*=0.
	\end{equation}
	\begin{enumerate}
		\item If, for each $x^*,s\in\mathbb X$ and $y^*,z^*,v\in\R^m$ satisfying
			$\dual{y^*}{v} \geq 0$, \eqref{eq:pseudo_coderivative_nonpolyhedral_x}, and
			\begin{equation}\label{eq:system_pseudo_coderivative_polyhedral_II}
				\begin{aligned}
						y^*
						&\in 
						\mathcal N_{\mathbf T(u)}(w_s(u,v))
						\cap
						\ker\nabla g(\bar x)^*,
						\\
						z^*
						&\in 
						\mathcal T_{
							\mathcal N_{\mathbf T(u)}(w_s(u,v))
							}
							(y^*),
				\end{aligned}
			\end{equation}
			there is some $\lambda\in\mathcal N_D(g(\bar x))$
			such that $x^*=\nabla g(\bar x)^*\lambda$,
			then $\Phi$ is asymptotically regular at $(\bar x,0)$ in direction $u$.
		\item If, for each $x^*,s\in\mathbb X$ and $y^*,z^*,v\in\R^m$ satisfying
			$\dual{y^*}{v} \geq 0$, \eqref{eq:pseudo_coderivative_nonpolyhedral_x}, and
			\eqref{eq:system_pseudo_coderivative_polyhedral_II}, there is some
			$\lambda\in\mathcal N_{\mathcal T_D(g(\bar x))}(\nabla g(\bar x)u)$
			such that $x^*=\nabla g(\bar x)^*\lambda$, then $\Phi$ is
			strongly asymptotically regular at $(\bar x,0)$ in direction $u$.			
	\end{enumerate}
\end{theorem}
\begin{proof}
	We start with the proof of the first assertion.
	With the aid of \cite[Theorem~3.2\,(d)]{BenkoMehlitz2022}, we easily see that
	\eqref{eq:CQ_pseudo_subregularity_polyhedral_II} yields
	$\ker D^*_2\Phi((\bar x,0);(u,0))=\{0\}$ in the present situation.
	Now, fix $x^*\in \widetilde D^*_2\Phi((\bar x,0);(u,0))(0)$.
	Then \cite[Theorem~3.2\,(c)]{BenkoMehlitz2022} shows the existence
	of $z^*\in\mathcal N_{\mathcal T_D(g(\bar x))}(\nabla g(\bar x)u)$ such that
	$x^*=\nabla g(\bar x)^*z^*$.
	In case where $x^*\in D^*_2\Phi((\bar x,0);(u,\alpha w))(\beta w)$ holds for 
	some $w\in\mathbb S_{\R^m}$ and $\alpha,\beta \geq 0$,
	\cite[Theorem~3.2\,(d)]{BenkoMehlitz2022} implies the existence of $s\in\mathbb X$ such that
	\eqref{eq:pseudo_coderivative_nonpolyhedral_x} and
	\eqref{eq:system_pseudo_coderivative_polyhedral_II} hold with $v:=\alpha w$
	and $y^*:=\beta w$,
	and this gives $\dual{y^*}{v}=\alpha\beta \norm{w}^2\geq 0$.
	Now, the postulated assumptions guarantee the existence of
	$\lambda\in\mathcal N_{\mathcal T_D(g(\bar x))}(\nabla g(\bar x)u)$ such that
	$x^*=\nabla g(\bar x)^*\lambda$.
	Respecting 
	\cref{lem:coderivatives_constraint_maps}\,\ref{item:constraint_maps_limiting_coderivative},
	this shows
	\eqref{eq:sufficient_condition_asym_reg_pseudo} with $\bar y:=0$ and $\gamma:=2$.
	Thus, \cref{thm:asymptotic_regularity_via_pseudo_coderivatives} yields
	that $\Phi$ is asymptotically regular at $(\bar x,0)$ in direction $u$.
	
	The second statement follows in analogous fashion while exploiting
	\cref{lem:coderivatives_constraint_maps}\,\ref{item:constraint_maps_directional_limiting_coderivative}.
\end{proof}

We note that the sufficient conditions for directional asymptotic
regularity stated in 
\cref{thm:asymptotic_regularity_pseudo_coderivative_non_polyhedral,thm:asymptotic_regularity_pseudo_coderivative_polyhedral_I,thm:asymptotic_regularity_pseudo_coderivative_polyhedral_II}
recover the constraint qualifications for M-stationarity we obtained 
in \cite[Section~4.2]{BenkoMehlitz2022} by a different approach.
In the remaining part of the paper, we prove that the assumptions of
\cref{thm:asymptotic_regularity_pseudo_coderivative_non_polyhedral}
are not stronger than FOSCMS
while the assumptions of
\cref{thm:asymptotic_regularity_pseudo_coderivative_polyhedral_I,thm:asymptotic_regularity_pseudo_coderivative_polyhedral_II}
are weaker than the so-called \emph{Second-Order Sufficient Condition for Metric Subregularity} (SOSCMS).

Given a point $\bar x\in\mathbb X$ with $(\bar x,0)\in\gph\Phi$,
\cref{lem:coderivatives_constraint_maps}\,\ref{item:constraint_maps_directional_limiting_coderivative} 
shows that the condition
\[
	u\in\mathbb S_{\mathbb X},\,
	\nabla g(\bar x)u\in\mathcal T_D(g(\bar x)),\,
	\nabla g(\bar x)^*y^*=0,\,
	y^*\in\mathcal N_D(g(\bar x);\nabla g(\bar x)u)
	\quad
	\Longrightarrow
	\quad
	y^*=0
\]
equals FOSCMS in the current setting.
In case where $D$ is locally polyhedral around $g(\bar x)$, the refined condition
\[
	\left.
	\begin{aligned}
	&u\in\mathbb S_{\mathbb X},\,
	\nabla g(\bar x)u\in\mathcal T_D(g(\bar x)),\,
	\nabla g(\bar x)^*y^*=0,\\
	&
	\nabla^2\innerprod{y^*}{g}(\bar x)[u,u]\geq 0,\,
	y^*\in\mathcal N_D(g(\bar x);\nabla g(\bar x)u)
	\end{aligned}
	\right\}
	\quad
	\Longrightarrow
	\quad
	y^*=0,
\]
is referred to as SOSCMS in the literature.
As these names suggest, both conditions are sufficient for metric subregularity of
$\Phi$ at $(\bar x,0)$, see \cite[Corollary~1]{GfrererKlatte2016}.
Particularly, they provide constraint qualifications for M-stationarity of local minimizers.
Moreover, validity of these conditions for a fixed direction $u\in\mathbb S_{\mathbb X}$,
denoted by FOSCMS$(u)$ and SOSCMS$(u)$, respectively,
is sufficient for metric subregularity of $\Phi$ at $(\bar x,0)$ in direction $u$.

We split our remaining considerations into the general and the polyhedral case.
\begin{proposition}\label{Pro:Milder_than_FOSCMS}
	Consider $(\xb,0) \in \gph \Phi$ and $u\in\mathbb S_{\mathbb X}$.
	Under FOSCMS$(u)$ all assumptions of \cref{thm:asymptotic_regularity_pseudo_coderivative_non_polyhedral} are satisfied.
\end{proposition}
\begin{proof}
Let $y^* \in \mathcal N_{D}(g(\xb);\nabla g(\xb) u)$
be such that $\nabla g(\xb)^* y^* = 0$.
Then FOSCMS$(u)$ yields $y^* = 0$ 
and so \eqref{eq:CQ_pseudo_subregularity_nonpolyhedral_II} is satisfied.
Moreover, we only need to show the remaining assertions for $y^* = 0$.

First, we claim that
\begin{equation}\label{eq:I_have_no_idea_how_to_name_this}
	D\mathcal N_{D}(g(\xb),0)(q) \ \Big( D_{\textrm{sub}}\mathcal N_{D}(g(\xb),0)(q)\Big)
	\ \subset \
	\mathcal N_D(g(\xb);q)
\end{equation}
holds for any $q \in \mathbb Y$ ($q \in \mathbb S_{\mathbb Y}$).
Indeed, let $\hat{z}^* \in D\mathcal N_{D}(g(\xb),0)(q)$.
The definition yields that there are sequences
$\{t_k\}_{k\in\N} \subset\R_+$ and $\{q_k\}_{k\in\N},\{\hat z_k^*\}_{k\in\N}\subset \mathbb Y$
with $t_k\searrow 0$, $q_k \to q$, $\hat{z}_k^* \to \hat{z}^*$,
and $t_k \hat{z}_k^* \in \mathcal N_D(g(\xb) + t_k q_k)$ for each $k\in\N$.
Since, for each $k\in\N$, $\mathcal N_D(g(\xb) + t_k q_k)$ is a cone, however, 
we get $\hat{z}_k^* \in \mathcal N_D(g(\xb) + t_k q_k)$,
and $\hat{z}^* \in \mathcal N_D(g(\xb);q)$ follows by robustness of the directional limiting
normal cone, see \cref{lem:robustness_directional_limiting_normals}.
The case $D_{\textrm{sub}}\mathcal N_{D}(g(\xb),0)(q)$ is almost identical.

Next, assume that $\nabla g(\xb)u \neq 0$ holds.
Suppose now that \eqref{eq:CQ_pseudo_subregularity_nonpolyhedral_Ib} is violated, i.e., 
there exists
$\hat z^* \in D_{\textrm{sub}}\mathcal N_{D}(g(\xb),0)(q)$ for $q:=\nabla g(\xb)u/\norm{\nabla g(\xb)u}$ with $\nabla g(\bar x)^*\hat z^*=0$.
By \eqref{eq:I_have_no_idea_how_to_name_this} and FOSCMS$(u)$ we thus get $\hat z^* = 0$
which is a contradition since $\hat z^* \in  \mathbb S_{\mathbb Y}$ by definition.
Similarly, in case $\nabla g(\bar x)u=0$, we can verify \eqref{eq:CQ_pseudo_subregularity_nonpolyhedral_Ia} which reduces to
\[
	\nabla g(\bar x)^*\hat z^*=0,\quad
	\hat{z}^* \in D\mathcal N_{D}(g(\xb),0)(0)
	\quad
	\Longrightarrow
	\quad
	\hat z^*=0.
\]
Applying \eqref{eq:I_have_no_idea_how_to_name_this} with $q:=0$, 
we get $\hat{z}^* \in \mathcal N_D(g(\xb))$ which again implies $\hat{z}^* = 0$
since FOSCMS$(u)$ corresponds to the Mordukhovich criterion 
due to $\nabla g(\xb)u = 0$.
Thus, we have shown that \eqref{eq:CQ_pseudo_subregularity_nonpolyhedral_Ia} or,
in case $\nabla g(\bar x)u\neq 0$, \eqref{eq:CQ_pseudo_subregularity_nonpolyhedral_Ib}
holds.

Validity of the last assumption follows immediately
since $z^* \in \mathcal N_D(g(\xb);\nabla g(\xb)u)$ follows from
\eqref{eq:I_have_no_idea_how_to_name_this}, and so we can just take $\lambda := z^*$
due to $y^*=0$.
\end{proof}

\begin{remark}\label{rem:forgetting_directional_information_in_CQ}
	Note that for $u\in\mathbb S_{\mathbb X}$ satisfying $\nabla g(\bar x)u\neq 0$, we have the trivial
	upper estimate
	$D_{\textup{sub}}\mathcal N_{D}(g(\xb),y^*)(\nabla g(\xb)u/\norm{\nabla g(\bar x)u}) 
	\subset D\mathcal N_{D}(g(\xb),y^*)(0)$,
	but keeping only this non-directional information, i.e.,
	relying only on \eqref{eq:CQ_pseudo_subregularity_nonpolyhedral_Ia},
	we cannot show that FOSCMS$(u)$ is sufficient for it to hold.
\end{remark}

\begin{proposition}\label{Pro:Milder_than_SOSCMS}
	Let $(\bar x,0)\in\gph\Phi$ as well as $u\in\mathbb S_{\mathbb X}$
	be fixed and assume that SOSCMS$(u)$ is valid.
	Furthermore, let $\mathbb Y:=\R^m$ and let $D$ be polyhedral
	locally around $g(\bar x)$.
	Then the following statements hold.
	\begin{enumerate}
		\item The assumptions of 
			\cref{thm:asymptotic_regularity_pseudo_coderivative_polyhedral_II} are satisfied.
		\item\label{case b} The assumptions of  
			\cref{thm:asymptotic_regularity_pseudo_coderivative_polyhedral_I} hold for $x^*\in\mathbb X$
			satisfying $\dual{x^*}{u}\geq 0$.
	\end{enumerate}
\end{proposition}
\begin{proof}
For the proof of the first statement, 
let $y^* \in \mathcal N_{\mathbf T(u)}(w_s(u,0))$ for some $s\in\mathbb X$
be such that $\nabla g(\xb)^* y^* = 0$.
We get
\begin{align*}
	\mathcal N_{\mathbf T(u)}(w_s(u,0)) 
	&=
	\mathcal N_{\mathcal T_D(g(\bar x))}(\nabla g(\bar x) u; w_s(u,0))
	\\
	&
	\subset 
	\mathcal N_{\mathcal T_D(g(\bar x))}(\nabla g(\bar x) u) \cap [w_s(u,0)]^{\perp}
	=
	\mathcal N_D(g(\bar x);\nabla g(\bar x)u)\cap[w_s(u,0)]^\perp
\end{align*}
due to (local) polyhedrality of $\mathcal T_D(g(\bar x))$ and $D$ 
from \cite[Lemma~2.1]{BenkoMehlitz2022}.
From $\nabla g(\xb)^* y^* = 0$ we thus obtain
\[
	\frac12 \nabla^2\langle y^*,g\rangle(\bar x)[u,u]
	=
	\frac12 \innerprod{\nabla^2 g(\xb)[u,u]}{ y^* }
	= 
	\innerprod{w_s(u,0)}{ y^* } = 0
\]
and SOSCMS$(u)$ yields $y^* = 0$ which shows validity of \eqref{eq:CQ_pseudo_subregularity_polyhedral_II}.

Next, for arbitrary $y^* \in \mathcal N_{\mathbf T(u)}(w_s(u,v))\cap\ker\nabla g(\bar x)^*$ with $s\in\mathbb X$ and 
$v\in\R^m$ satisfying $\innerprod{y^*}{v}\geq 0$, we get
$y^*\in\mathcal N_D(g(\bar x);\nabla g(\bar x)u)\cap[w_s(u,v)]^\perp$ and
\[
	\frac12 \nabla^2\langle y^*,g\rangle(\bar x)[u,u]
	= 
	\innerprod{w_s(u,v)}{ y^*}+ \innerprod{v}{y^*}
	=
	 \innerprod{v}{ y^* }
	\geq 
	0,
\]
so SOSCMS$(u)$ can still be applied to give $y^* = 0$.
Now, for $z^* \in \mathcal T_{\mathcal N_{\mathbf T(u)}(w_s(u,v))}(0)$, 
we get 
$z^* \in \mathcal N_{\mathbf T(u)}(w_s(u,v)) \subset \mathcal N_D(g(\bar x);\nabla g(\bar x)u)$
and so under SOSCMS$(u)$ we can always take $\lambda := z^*$ since $y^* = 0$.

In order to prove the second claim, let $y^*,z^*\in\R^m$
satisfy the requirements of condition \eqref{eq:CQ_pseudo_subregularity_polyhedral_I}.
By (local) polyhedrality of $\mathcal N_{\mathcal T_D(g(\bar x))}(\nabla g(\bar x) u)$
and $D$, for sufficiently small $\alpha > 0$, we obtain
\[
	y^* + \alpha z^* 
	\in 
	\mathcal N_{\mathcal T_D(g(\bar x))}(\nabla g(\bar x) u) 
	\subset 
	\mathcal N_D(g(\bar x)) \cap [\nabla g(\bar x) u]^{\perp},
\]
see \cite[Lemma~2.1]{BenkoMehlitz2022} again.
Due to $\nabla g(\bar x)^* y^* = 0$, we also get 
\[
	0
	=
	\innerprod{u}{\nabla g(\bar x)^*(y^*+\alpha z^*)}
	=
	\alpha\,\innerprod{u}{\nabla g(\bar x)^* z^*},
\]
i.e., $\innerprod{u}{\nabla g(\bar x)^* z^* }= 0$.
Taking into account the equation $\nabla^2\langle y^*,g\rangle(\bar x)(u) + \nabla g(\xb)^* z^* = 0$, we get
\[
	\nabla^2\langle y^*,g\rangle(\bar x)[u,u]
	=
	\innerprod{u}{\nabla^2\langle y^*,g\rangle(\bar x)(u)}
	= 
	-\innerprod{u}{\nabla g(\bar x)^* z^*}
	= 
	0
\]
and SOSCMS$(u)$ yields $y^* = 0$ since $y^*\in\mathcal N_D(g(\bar x);\nabla g(\bar x)u)$.

Finally, assume that there are
$x^*\in\mathbb X$, $y^*\in\mathcal N_{\mathcal T_D(g(\bar x))}(\nabla g(\bar x)u)\cap\ker\nabla g(\bar x)^*$, and
$z^*\in\mathcal T_{\mathcal N_{\mathcal T_D(g(\bar x))}(\nabla g(\bar x)u)}(y^*)$  such that $\dual{x^*}{u}\geq 0$ and
$x^*=\nabla^2\langle y^*,g\rangle(\bar x)(u)+\nabla g(\bar x)^* z^*$.
Exploiting the above arguments, we find
\[
	\nabla^2\langle y^*,g\rangle(\bar x)[u,u]
	=
	\innerprod{u}{\nabla^2\langle y^*,g\rangle(\bar x)(u)}
	=
	-\innerprod{u}{\nabla g(\bar x)^* z^*}+\innerprod{u}{x^*}
	=
	\innerprod{u}{x^*}
	\geq 
	0
\]
and due to $y^*\in\mathcal N_D(g(\bar x);\nabla g(\bar x)u)$, SOSCMS$(u)$ yields $y^*=0$.
Thus, we also have $z^*\in\mathcal N_D(g(\bar x);\nabla g(\bar x)u)$ and we can choose $\lambda:=z^*$
again.
\end{proof}

We immediately arrive at the following corollary.

\begin{corollary}
	The constraint mapping $\Phi$ is strongly asymptotically regular at $(\xb,0) \in \gph \Phi$ in direction $u\in\mathbb S_{\mathbb X}$ if
	FOSCMS$(u)$ holds or if $\mathbb Y:=\R^m$, $D$ is locally polyhedral around $g(\xb)$, and SOSCMS$(u)$ holds.
\end{corollary}

\begin{remark}
	Note that showing that SOSCMS implies strong directional asymptotic regularity was only possible via the refined
	conditions from \cref{thm:asymptotic_regularity_pseudo_coderivative_polyhedral_II}.
	In the context of M-stationarity, however, also the simpler, more restrictive conditions from
	\cref{thm:asymptotic_regularity_pseudo_coderivative_polyhedral_I} can be useful.
	Indeed, let $\bar x\in\mathcal F$ be a local minimizer of \eqref{eq:nonsmooth_problem}.
	In order to justify M-stationarity of $\bar x$, it is sufficient to verify the assumptions from
	\cref{thm:asymptotic_regularity_pseudo_coderivative_polyhedral_I} for $x^*\in - \partial \varphi(\bar x)$
	and, taking into account \cite[Remark~4.7]{BenkoMehlitz2022},
	we only need to consider
	\[
		-x^*
		\in
		\partial \varphi(\bar x;(u,\mu))
		:=
		\{x^*\in\mathbb X \,|\, (x^*,-1) \in \mathcal N_{\epi \varphi}((\bar x,\varphi(\bar x));(u,\mu))\}
	\]
	for some $\mu \leq 0$, where $\partial \varphi(\bar x;(u,\mu))$
	denotes the (geometric) subdifferential of $\varphi$ at $\bar x$ in direction $(u,\mu)$, see \cite{BenkoGfrererOutrata2019}.
	Then, whenever $\epi\varphi$ is so-called semismooth* at $(\bar x,\varphi(\bar x))$, we directly get $\innerprod{x^*}{u}= -\mu \geq 0$
	and \cref{Pro:Milder_than_SOSCMS}~\ref{case b} can be applied.

	Note that, as recently introduced in \cite{GfrererOutrata2021}, a closed
	set $Q\subset\mathbb X$ is called
	\emph{semismooth*} at $x\in Q$ if for all $w\in\mathcal T_Q(x)$ and
	$\eta\in\mathcal N_Q(x;w)$, we have $\innerprod{\eta}{w}=0$.
	Let us briefly mention that broad classes of sets enjoy the semismoothness*
	property, see \cite[Section~3]{GfrererOutrata2021}. 
	Exemplary, each union of finitely many closed, convex sets and graphs, epigraphs, as well as
	hypographs of continuously differentiable mappings are semismooth* everywhere.
\end{remark}

The following example shows that our new conditions are in fact strictly milder than SOSCMS.
We conjecture that they are also strictly milder than FOSCMS but, unfortunately, have no example
available which shows this.

\begin{example}\label{ex:milder_than_SOSCMS}
	Let $g\colon\R \to \R^2$ and $D \subset \R^2$ be given by
	$g(x) := (x,-x^2)$ for all $x\in\R$ and $D:=(\R_+ \times \R) \cup (\R \times \R_+)$.
	Observe that $D$ is a polyhedral set.
	We consider the constraint map $\Phi\colon\R\tto\R^2$ given by $\Phi(x):=g(x)-D$ for all $x\in\R$.
	We note that $\Phi^{-1}(0)=[0,\infty)$ holds.
	Hence, fixing $\bar x:=0$, we can easily check that $\Phi$ is metrically
	subregular at $(\bar x,0)$ in direction $1$ but not in direction $-1$.
	Hence, FOSCMS and SOSCMS must be violated.
	
	First, we claim that all the assumptions from \cref{thm:asymptotic_regularity_pseudo_coderivative_polyhedral_I}
	(and, hence, also \cref{thm:asymptotic_regularity_pseudo_coderivative_polyhedral_II}) are satisfied for $u=\pm 1$.
	Thus, let us fix $u=\pm 1$, $y^*\in\mathcal N_{\mathcal T_D(g(\xb))}(\nabla g(\bar x)u)\cap\ker\nabla g(\bar x)^*$,
	and $z^*\in\mathcal T_{\mathcal N_{\mathcal T_D(g(\xb))}(\nabla g(\bar x)u)}(y^*)$ such that 
	$\nabla^2\langle y^*,g\rangle(\bar x)(u)+\nabla g(\bar x)^* z^*=0$.
	From $y^*\in\ker\nabla g(\bar x)^*$ we have $y_1^*=0$.
	Furthermore, we have $\nabla g(\bar x)u=(u,0)$, 
	$\nabla^2\langle y^*,g\rangle(\bar x)(u)=-2y_2^*u$, and
	\[
		\mathcal N_{\mathcal T_D(g(\xb))}(\nabla g(\bar x)u) 
		=
		\begin{cases} 
		\{0\} \times \R_-	& u=-1,\\
		\{(0,0)\}			& u=1.
		\end{cases}
	\]
	Thus, for $u=1$, condition \eqref{eq:CQ_pseudo_subregularity_polyhedral_I} holds trivially. 
	For $u=-1$, we fix $y^*_2\leq 0$ and, thus,
	$\mathcal T_{\mathcal N_{\mathcal T_D(g(\xb))}(\nabla g(\bar x)u)}(y^*) \subset \{0\} \times \R$, i.e., $z_1^*=0$.
	Thus, from $-2 y_2^*u + z_1^* = 0,$
	we deduce $y_2^*=0$, and \eqref{eq:CQ_pseudo_subregularity_polyhedral_I} follows.
	In order to check the second assumption of \cref{thm:asymptotic_regularity_pseudo_coderivative_polyhedral_I},
	we fix $x^*\in\R$, $u=\pm 1$, $y^*\in\mathcal N_{\mathcal T_D(g(\xb))}(\nabla g(\bar x)u)\cap\ker\nabla g(\bar x)^*$,
	and $z^*\in\mathcal T_{\mathcal N_{\mathcal T_D(g(\xb))}(\nabla g(\bar x)u)}(y^*)$ such that 
	$x^* = \nabla^2\langle y^*,g\rangle(\bar x)(u)+\nabla g(\bar x)^* z^*$.
	In case $u=1$, we have $y^*=0$ from above. This yields $z^*\in\mathcal N_{\mathcal T_D(g(\xb))}(\nabla g(\bar x)u)$,
	and we can choose $\lambda:=z^*$ to find $x^* = \nabla g(\bar x)^*\lambda$ as well as $\lambda\in\mathcal N_{\mathcal T_D(g(\xb))}(\nabla g(\bar x)u)$.
	Thus, let us consider $u=-1$. Then we find $y_2^*\leq 0$ and $z_1^*=0$ from above.
	Next from $x^* = -2 y_2^*u + z_1^* = 2 y_2^* \leq 0$ we can choose $\lambda := (x^*,0) \in \mathcal N_D(g(\bar x))$
	to get $\nabla g(\bar x)^*\lambda=x^*$.

	Note, however, that $\lambda = (x^*,0) \notin \mathcal N_{\mathcal T_D(g(\xb))}(\nabla g(\bar x)u) = \{0\} \times \R_-$
	unless $x^* = 0$.
	
	Regarding the assumptions of \cref{thm:asymptotic_regularity_pseudo_coderivative_non_polyhedral},
	let us just mention, without providing all the details, that \eqref{eq:CQ_pseudo_subregularity_nonpolyhedral_Ia} and \eqref{eq:CQ_pseudo_subregularity_nonpolyhedral_Ib}
	fail since the graphical (sub-) derivative is too large.
	Particularly, this clarifies that the first assumption is not necessary e.g.\ in the polyhedral setting,
	but not because it would be automatically satisfied.

\if{
	Let us first show that the first assumption of \cref{Pro:M-stat_via_second_order}\,\ref{item:pseudo_stationarity_constraint_maps}
	fails. Therefore, let us first mention that
	\[
		\gph\mathcal N_D
		=
		(\R_-\times\{0\})\times(\{0\}\times\R_-)
		\cup
		(\{0\}\times\R_-)\times(\R_-\times\{0\})
		\cup
		D\times(\{0\}\times\{0\})
	\]
	is valid. For $u:=-1$, we fix $y^*\in\mathcal N_D(g(\bar x);\nabla g(\bar x)u)\cap\ker\nabla g(\bar x)^*=\{0\}\times\R_-$
	and $\hat z^*\in D_\textup{sub}\mathcal N_D(g(\bar x),y^*)(\nabla g(\bar x)u)\cap\ker \nabla g(\bar x)^*$.
	Then we directly have $y_1^*=\hat z_1^*=0$ and $y_2^*\leq 0$
	In case $u>0$, we find $y^*=0$ and the first component of $\nabla g(\bar x)u$ is positive.
	Let us consider the case $y_2^*<0$.
	Then, we find that $\gph \mathcal N_D$ coincides with
	$(\R_-\times\{0\})\times(\{0\}\times\R_-)$ locally around $(g(\bar x),y^*)$.
	Due to \cref{rem:graph_cartesian_product}, this yields 
	$\gph D_\textup{sub}\mathcal N_D(g(\bar x),y^*)=(\R_-\times\{0\})\times(\{0\}\times\R)$.
	Particularly, we find $\hat z_1^*=0$ which we already know. However, we cannot deduce $\hat z_2^*=0$ so the first assumption from
	\cref{Pro:M-stat_via_second_order}\,\ref{item:pseudo_stationarity_constraint_maps} fails.
	One can check, however, that the second and third assumption are valid.
}\fi
\end{example}

\section{Concluding remarks}\label{sec:conclusions}

In this paper, we introduced directional notions of asymptotic regularity for set-valued
mappings. These conditions have been shown to serve as constraint qualifications guaranteeing
M-stationarity of local minimizers in nonsmooth optimization. These new qualification
conditions have been embedded into the landscape of constraint qualifications which are
already known from the literature, and we came up with the impression that these conditions
are comparatively mild.
Noting that directional asymptotic regularity might be difficult to check in practice,
we then focused on the derivation of applicable sufficient conditions for its validity.
First, we suggested directional notions of pseudo- and quasi-normality for that purpose
which have been shown to generalize related concepts for geometric constraint systems to
arbitrary set-valued mappings. 
Second, with the aid of so-called super- and pseudo-coderivatives, sufficient conditions
for the presence of directional asymptotic regularity for geometric constraint systems
in terms of first- and second-order derivatives of the associated mapping as well as
standard variational objects associated with the underlying set were derived. 
These sufficient conditions turned out to recover some of 
our findings from \cite{BenkoMehlitz2022},
and we showed that they are not stronger than FOSCMS and SOSCMS.
In this paper, we completely neglected to study the potential value of directional
asymptotic regularity in numerical optimization 
which might be a promising topic of future research.
Furthermore, it has been shown in \cite{Mehlitz2020a} that non-directional asymptotic
regularity can be applied nicely as a qualification condition in the limiting variational
calculus. Most likely, directional asymptotic regularity may play a similar role 
in the directional limiting calculus. 

\subsection*{Acknowledgements}
 	The research of the first author was supported by the Austrian Science Fund (FWF) under grant P32832-N.


\end{document}